\documentclass{amsart}

\newtheorem{theorem}{Theorem}[section]
\newtheorem{lemma}[theorem]{Lemma}

\newtheorem{proposition}[theorem]{Proposition}

\theoremstyle{definition}
\newtheorem{definition}[theorem]{Definition}

\theoremstyle{remark}
\newtheorem{remark}[theorem]{Remark}

\usepackage{graphicx}
\usepackage{color}
\usepackage{amsmath}
\usepackage{amsfonts}
\usepackage{amssymb}
\newcommand{\be}{\begin{equation}}
\newcommand{\ee}{\end{equation}}

\newcommand{\Om}{\Omega}

\newcommand{\lc}{{\stackrel{\scriptscriptstyle{LC}}{\Gamma}}}

\newcommand{\dz}{\wedge}

\newcommand{\ba}{\begin{array}}

\newcommand{\ea}{\end{array}}

\newcommand{\beq}{\begin{eqnarray}}

\newcommand{\eeq}{\end{eqnarray}}

\newtheorem{lm}{lemma}

\newtheorem{thee}{theorem}

\newtheorem{proo}{proposition}

\newtheorem{co}{corollary}

\newtheorem{rem}{remark}

\newtheorem{deff}{definition}

\newcommand{\bd}{\begin{deff}}

\newcommand{\ed}{\end{deff}}

\newcommand{\bl}{\begin{lm}}

\newcommand{\el}{\end{lm}}

\newcommand{\bp}{\begin{proo}}

\newcommand{\ep}{\end{proo}}

\newcommand{\bt}{\begin{thee}}

\newcommand{\et}{\end{thee}}

\newcommand{\bc}{\begin{co}}

\newcommand{\ec}{\end{co}}

\newcommand{\brm}{\begin{rem}}

\newcommand{\erm}{\end{rem}}

\newcommand{\der}{{\rm d}}

\newcommand{\sgn}{\mathrm{sgn}}

\usepackage{t1enc}
\def\frak{\mathfrak}

\newcommand{\newc}{\newcommand}

\renewcommand{\exp}{\operatorname{exp}}

\let\ccdot\cdot
\def\cdot{\hbox to 2.5pt{\hss$\ccdot$\hss}}

\newc{\aR}{\mbox{\boldmath{$ R$}}}
\newc{\aS}{\mbox{\boldmath{$ S$}}}
\newc{\aT}{\mbox{\boldmath{$ T$}}}
\newc{\aW}{\mbox{\boldmath{$ W$}}}

\newc{\aK}{\mbox{\boldmath{$ K$}}}
\newc{\aL}{\mbox{\boldmath{$ L$}}}
\newcommand{\bbC}{\mathbb{C}}
\newcommand{\bbH}{\mathbb{H}}
\newcommand{\bbO}{\mathbb{O}}
\newcommand{\bbK}{\mathbb{K}}

\usepackage{amssymb}
\usepackage{amscd}

\newc{\obstrn}[2]{B^{#1}_{#2}}

\newcommand{\rpl}                         % +) or <+
{\mbox{$
\begin{picture}(12.7,8)(-.5,-1)
\put(0,0.2){$+$}
\put(4.2,2.8){\oval(8,8)[r]}
\end{picture}$}}

\newcommand{\lpl}                         % (+ or +>
{\mbox{$
\begin{picture}(12.7,8)(-.5,-1)
\put(2,0.2){$+$}
\put(6.2,2.8){\oval(8,8)[l]}
\end{picture}$}}

\usepackage{ifthen}
\newcommand{\bbM}{\mathbb{M}}
\newcommand{\bbR}{\mathbb{R}}
\newcommand{\sog}{\mathbf{SO}}

\newcommand{\slg}{\mathbf{SL}}
\newcommand{\og}{\mathbf{O}}
\newcommand{\soa}{\frak{so}}

\newcommand{\sua}{\frak{su}}
\newcommand{\spa}{\frak{sp}}

\newcommand{\spg}{\mathbf{Sp}}
\newcommand{\sug}{\mathbf{SU}}
\newcommand{\ug}{\mathbf{U}}

\newcommand{\thet}{\tilde{\theta}}

\newc{\tensor}[1]{#1}
\newc{\Mvariable}[1]{\mbox{#1}}
\newc{\down}[1]{{}_{#1}}
\newc{\up}[1]{{}^{#1}}

\newc{\JulyStrut}{\rule{0mm}{6mm}}
\newc{\midtenPan}{\mbox{\sf S}}
\newc{\midten}{\mbox{\sf T}}
\newc{\midtenEi}{\mbox{\sf U}}
\newc{\ATen}{\mbox{\sf E}}
\newc{\BTen}{\mbox{\sf F}}
\newc{\CTen}{\mbox{\sf G}}

\def\sideremark#1{\ifvmode\leavevmode\fi\vadjust{\vbox to0pt{\vss% the remark
 \hbox to 0pt{\hskip\hsize\hskip1em%                          will appear only
 \vbox{\hsize3cm\tiny\raggedright\pretolerance10000%          on the side
 \noindent #1\hfill}\hss}\vbox to8pt{\vfil}\vss}}}%
\newcommand{\bgw}{{\textstyle \bigwedge}}
\newcommand{\bgs}{{\textstyle \bigodot}}
\newcommand{\bgt}{{\textstyle \bigotimes}}

\numberwithin{equation}{section}

\newcounter{romenumi}
\newcommand{\labelromenumi}{(\roman{romenumi})}

\newcommand{\ten}{\Upsilon}

\begin{document}
\title{Distinguished dimensions for special Riemannian geometries}

\author{Pawe\l~ Nurowski} \address{Instytut Fizyki Teoretycznej,
Uniwersytet Warszawski, ul. Hoza 69, Warszawa, Poland}
\email{nurowski@fuw.edu.pl} \thanks{This research was supported by
the KBN grant 1 P03B 07529}

\date{\today}

\begin{abstract} The paper is based on relations between a
ternary symmetric form defining the $\sog(3)$ geometry in dimension
five and Cartan's works on isoparametric hypersurfaces in spheres. As observed by Bryant such a ternary form exists only in dimensions $n_k=3k+2$,
where $k=1,2,4,8$. In these dimensions it reduces the orthogonal
group to the subgroups $H_k\subset \sog(n_k)$, with $H_1=\sog(3)$,
$H_2=\sug(3)$, $H_4=\spg(3)$ and $H_8={\bf F_4}$. This enables
studies of special Riemannian geometries with structure groups $H_k$
in dimensions $n_k$.\\
The necessary and sufficient conditions for the $H_k$ geometries to
admit the characteristic connection are given. As an illustration
nontrivial examples of $\sug(3)$ geometries in dimension 8 admitting
characteristic connection are provided. Among them are the examples
having nonvanishing torsion and satisfying Einstein equations with
respect to either the Levi-Civita or the characteristic
connections.\\
The torsionless models for the $H_k$ geometries have the respective
symmetry groups $G_1=\sug(3)$,
$G_2=\sug(3)\times\sug(3)$, $G_3=\sug(6)$ and
$G_4={\bf E}_6$. The groups $H_k$ and $G_k$
constitute a part of the `magic square' for Lie groups. The `magic
square' Lie groups suggest studies of ten other classes of special
Riemannian geometries. Apart from the two exceptional cases, they
have the structure groups $\ug(3)$, ${\bf S}(\ug(3)\times\ug(3))$,
$\ug(6)$, ${\bf E}_6\times\sog(2)$, $\spg(3)\times\sug(2)$,
$\sug(6)\times\sug(2)$, $\sog(12)\times\sug(2)$ and ${\bf
E}_7\times\sug(2)$ and should be considered in respective dimensions
12, 18, 30, 54, 28, 40, 64 and 112. The two `exceptional' cases are:
$\sug(2)\times\sug(2)$ geometries in dimension 8 and
$\sog(10)\times\sog(2)$ geometries in
dimension 32.\\
The case of $\sug(2)\times\sug(2)$ geometry in dimension 8 is
examined closer. We determine the tensor that reduces $\sog(8)$ to
$\sug(2)\times\sug(2)$ leaving the more detailed studies of the
geometries based on the magic square ideas to the forthcoming paper.
%Some homogeneous examples of such geometries in dimension 8 and
%... are given.

\vskip5pt\centerline{\small\textbf{MSC classification}: 53A40, 53B15,
  53C10}\vskip15pt
\end{abstract}
\maketitle
%*************
\tableofcontents
%*************
\section{Introduction}
In a recent paper \cite{bobi} we
studied 5-dimensional manifolds $(M^5,g,\ten)$ equipped with the
Riemannian metric tensor
$g_{ij}$ and a 3-tensor $\ten_{ijk}$ such that
\begin{itemize}
\item[i)] %it was totally symmetric:
$\ten_{ijk}=\ten_{(ijk)}$,
\vskip1mm
\item[ii)]%it was trace free:
$\ten_{ijj}=0$,
\vskip1mm
\item[iii)] $\ten_{jki}\ten_{lmi}+\ten_{lji}\ten_{kmi}+\ten_{kli}\ten_{jmi}=g_{jk}g_{lm}+g_{lj}g_{km}+g_{kl}g_{jm}.$
%it satisfied the following condition:
\end{itemize}
It turns out that the quadratic condition iii) selects from all the
{\it symmetric totally trace free} 3-tensors in $\bbR^5$ only 
such whose stabilizer in $\sog(5)$ is the {\it irreducible} $\sog(3)$.

The geometry of Riemannian 5-manifolds $(M^5,g,\ten)$ is
particularly interesting if the tensor $\ten$ satisfies the {\it
nearly integrability} \cite{bobi} condition
\be
(\stackrel{\scriptscriptstyle{LC}}{\nabla}_v \ten)(v,v,v)\equiv 0.\label{nip}
\ee
In such case $(M^5,g,\ten)$ is naturally equipped with a unique metric
$\soa(3)$-valued connection $\nabla^T$ whose torsion $T$ is {\it totally skew
  symmetric}. We call this connection the {\it characteristic}
connection of the nearly integrable geometry $(M^5,g,\ten)$.
Existence and uniqueness of the characteristic connection enables a
classification of the nearly integrable geometries $(M^5,g,\ten)$
according to the algebraic properties of $T$ and the curvature $K$
of $\nabla^T$. In Ref. \cite{bobi} examples were given of the nearly
integrable geometries $(M^5,g,\ten)$ with the characteristic
connection $\nabla^T$ having $K$ and $T$ of all the possible types
from the above mentioned classification. In particular, a
7-parameter family of such geometries admitting at each point two
$\sog(3)$-invariant vector spaces of $\nabla^T$-covariantly constant
spinors were given. However, in this family of examples, the
characteristic $\nabla^T$ connection was flat.

Properties of $\ten$ resemble a bit properties of the tensor $J$,
defining an almost hermitian structure on a Riemannian manifold. For
example, condition iii) for $\ten$ is an algebraic condition of the
same sort as the almost hermitian condition \be
J_{ij}J_{jk}=-g_{ik}\label{j} \ee for $J$. Also, the nearly
integrable condition (\ref{nip}) for $\ten$ is similar to the nearly
K\"ahler condition
$$(\stackrel{\scriptscriptstyle{LC}}{\nabla}_v J)(v)\equiv 0$$ for
$J$. Since the almost hermitian condition (\ref{j}) imposes a severe
restriction on the dimension of the manifold to be {\it even}, a
natural question arises if there are some restrictions on the
dimensions of $\bbR^n$ in which one can have a tensor with the
properties i)-iii). More precisely we ask the following
question:\\

\noindent
\centerline{\it In which dimensions $n$ the Euclidean space $(\bbR^n,g)$}
\centerline{\it can be equipped 
with a tensor $\ten$ satisfying conditions i)-iii)?}\\

\noindent
It is rather easy to show that dimensions $n\leq 4$ do not admit such
tensor. Following \cite{bobi} we know that in dimension $n=5$ the
tensor $\ten$ may be defined by
\be
\ten_{ijk}x^ix^jx^k=\tfrac{1}{2}\det
\begin{pmatrix}
x^5-\sqrt{3}x^4&\sqrt{3}x^3&\sqrt{3}x^2\\
\sqrt{3}x^3&x^5+\sqrt{3}x^4&\sqrt{3}x^1\\
\sqrt{3}x^2&\sqrt{3}x^1&-2x^5
\end{pmatrix}.\label{ten}
\ee
Thus, $\ten$ is defined as a tensor whose
components are coefficients of the homogeneous polynomial of third
degree obtained as the
determinant of a generic $3\times 3$ real symmetric trace free
matrix.

\section{Dimensions 5, 8, 14 and 26}\label{rreppee}

Robert Bryant \cite{bryant} remarks that other dimensions in which
$\ten$ with properties i)-iii) surely exists are: $n=8$,
$n=14$ and $n=26$. This is essentially due to the fact that numbers 5,8,14 and
26 are values of the sequence $n_k=3 k+2$ for $k=1,2,4$ and $8$,
respectively. These four values of $k$ corresponds to the only
possible dimensions of the normed division algebras $\bbR$, $\bbC$,
$\bbH$ and $\bbO$. To fully explain Bryant's remark we need some
preparations.

Let $\bbK=\bbR,\bbC,\bbH$ or $\bbO$ and let $A\in M_{3\times 3}(\bbK)$
be a {\it hermitian} $3\times 3$ matrix with entries in $\bbK$. The
word `hermitian' here means that
the entries $a_{ij}$ and $a_{ji}$ of $A$ are
mutually conjugate in $\bbK$, i.e. $a_{ji}=\overline{a}_{ij}$ for
$i,j=1,2,3$. In particular, the entries $a_{11},a_{22},a_{33}\in\bbR$.

We may formally write
$$\det A=\sum_{\pi\in S_3} \sgn\pi~
a_{1\pi(1)}a_{2\pi(2)}a_{3\pi(3)},$$
which after expansion reads:
$$
\det A=a_{11}a_{22}a_{33}-a_{12}a_{21}a_{33}-a_{13}a_{22}a_{31}
-a_{11}a_{23}a_{32}+a_{13}a_{21}a_{32}+a_{12}a_{23}a_{31}.
$$
Note that despite of the possible {\it noncommutativity}, or even
{\it nonassociativity}, of the product, the values of the first {\it four}
monomials in the above formal expression are well defined. This is because
among the three factors in each of the four monomials, at least one is a real
number $a_{ii}$, the other two being either
both real (in the first term) or conjugate to each other
(in the remaining three terms). Thus, the values
of these four monomials are real numbers and do not depend on the
order of their factors and the order of the multiplication. Passing to
the last two terms in the formula for $\det A$ we see, that {\it
  a'priori} there is a huge number of possibilities to order the
factors and the brackets in these two terms. But the requirement that the
sum of these terms is {\it real} reduces this huge number to only 12
possibilities. It turns out that out of these 12 possibilities only
{\it two} are really different. They are all equal either to
$(a_{12}a_{23})a_{31}+a_{13}(a_{32}a_{21})$ or to
$(a_{21}a_{32})a_{13}+a_{31}(a_{23}a_{12})$. Note that the first
expression becomes the second under the transformation
$A\to\overline{A}$. Moreover, such transformation does not affect
the values of the first four terms in $\det A$. Summing
up we have the following Lemma.

\begin{lemma}
Given a hermitian matrix $A\in M_{3\times 3}(\bbK)$ with entries $a_{ij}\in\bbK$ such that
$a_{ji}=\overline{a}_{ij}$, $i,j=1,2,3$, where
$\bbK=\bbR,\bbC,\bbH,\bbO$, there are only two possibilities to assign
a real value to the Weierstrass formula
$$\det A=\sum_{\pi\in S_3} \sgn\pi~
a_{1\pi(1)}a_{2\pi(2)}a_{3\pi(3)}$$ for the determinant of $A$.
These two possible values are given by
$${\det}_1 A=a_{11}a_{22}a_{33}-a_{12}a_{21}a_{33}-a_{13}a_{22}a_{31}
-a_{11}a_{23}a_{32}+a_{13}(a_{32}a_{21})+(a_{12}a_{23})a_{31}$$
or by
$${\det}_2 A=a_{11}a_{22}a_{33}-a_{12}a_{21}a_{33}-a_{13}a_{22}a_{31}
-a_{11}a_{23}a_{32}+a_{31}(a_{23}a_{12})+(a_{21}a_{32})a_{13}.$$
In general $\det_1 A\neq\det_2 A$ if $\bbK=\bbH$ or $\bbO$, but
$\det_1 A\to\det_2 A$ when $A\to\overline{A}$.
\end{lemma}

Let $(e_0,e_1,e_2,e_3,e_4,e_5,e_6,e_7)$ be the unit octonions. We have: $e_0^2=1=-e_\mu^2$, $e_\mu e_{\mu+1}=e_{\mu+3}$, $e_\mu
e_\nu=-e_\nu e_\mu$, $\mu\neq\nu=1,2,3,...,7$, with additional
relations resulting from the cyclic permutation of each triple
$(e_\mu,e_{\mu+1},e_{\mu+3})$.

It is convenient to introduce
\begin{eqnarray*}
X^1&=&x^1 e_0+x^6 e_1+x^9 e_2+x^{10}e_4+x^{15} e_3+x^{16} e_5+x^{17}
e_6+x^{18}e_7,\\
X^2&=&x^2 e_0+x^7 e_1+x^{11} e_2+x^{12}e_4+x^{19}
e_3+x^{20} e_5+x^{21} e_6+x^{22}e_7,\\
X^3&=&x^3 e_0+x^8 e_1+x^{13} e_2+x^{14}e_4+x^{23} e_3+x^{24}
e_5+x^{25} e_6+x^{26}e_7.
\end{eqnarray*}
Then $X^1, X^2, X^3$ are three generic octonions. We can consider them
to be the generic quaternions if $x^I=0$ for all
$I=15,16,,17,...,26$, and three generic complex numbers if $x^I=0$ for all
$I=9,10,,11,...,26$. If $x^I=0$ for all
$I=6,7,,8,...,26$ then $X^1, X^2, X^3$ are three generic real numbers.
Using this we define a $3\times 3$ hermitian matrix
\be
A=\begin{pmatrix}
x^5-\sqrt{3}~x^4&\sqrt{3}~X^3&\sqrt{3}~X^2\\
\sqrt{3}~\overline{X}^3&x^5+\sqrt{3}~x^4&\sqrt{3}~X^1\\
\sqrt{3}~\overline{X}^2&\sqrt{3}~\overline{X}^1&-2x^5
\end{pmatrix}\label{maca}
\ee
in full analogy to the matrix entering the formula (\ref{ten}). Now,
 we have two `characteristic polynomials': $\det_1 (A-\lambda I)$ and
 $\det_2 (A-\lambda I)$. They can be written as:
\begin{eqnarray*}
{\det}_1(A-\lambda
 I)&=&-\lambda^3-3g(\vec{X},\vec{X})\lambda+2\ten_1(\vec{X},\vec{X},\vec{X}),\\
{\det}_2(A-\lambda
 I)&=&-\lambda^3-3g(\vec{X},\vec{X})\lambda+2\ten_2(\vec{X},\vec{X},\vec{X}),
\end{eqnarray*}
where $(\vec{X})^I=x^I$, $I=1,2,3,..,n_k=3k+2$ and $k=1,2,4,8$.
%\begin{eqnarray*}
%I&=&1,2,3,...,5 ~\,\,{\rm iff}~ k=1,~\bbK=\bbR,\\
%I&=&1,2,3,...,8 ~\,\,{\rm iff}~ k=2,~\bbK=\bbC,\\
%I&=&1,2,3,...,14 ~{\rm iff}~ k=3,~\bbK=\bbH,\\
%I&=&1,2,3,...,26 ~{\rm iff}~ k=4,~\bbK=\bbO.
%\end{eqnarray*}
The bilinear form is:
$$g(\vec{X},\vec{X})=(x^1)^2+
(x^2)^2+...+(x^{n_k})^2=g_{IJ}x^Ix^J$$
and the two ternary forms are
\[\ten_1(\vec{X},\vec{X},\vec{X})=\tfrac{1}{2}{\rm det}_1A\quad\quad {\rm and}\quad\quad
\ten_2(\vec{X},\vec{X},\vec{X})=\tfrac{1}{2}{\rm det}_2A.
\]
Now, we have the following proposition, which is our formulation of
Bryant's \cite{bryant} remark.
\begin{proposition}
If $I,J,K=1,2,3,...,n_k=3k+2$, $k=1,2,4,8$, then the tensors
$\ten_{IJK}^1$ and $\ten_{IJK}^2$ given by
$$\ten_{IJK}^a=\tfrac{1}{6}\frac{\partial^3\ten_a(\vec{X},\vec{X},\vec{X})}{\partial
  x^I\partial x^J\partial
  x^K}\quad\quad\quad a=1,2,
%^2\ten_{IJK}=\tfrac{1}{6}\frac{\partial^3\ten_2(\vec{X},\vec{X},\vec{X})}{\partial  x^I\partial x^J\partial x^K}
$$
satisfy
\begin{itemize}
\item[i)] %it was totally symmetric:
$\ten_{IJK}^a=\ten_{(IJK)}^a$,
\item[ii)]
%it was trace free:
$\ten_{IJJ}^a=0$,
\item[iii)]
$\ten_{JKI}^a\ten_{LMI}^a+\ten_{LJI}^a\ten_{KMI}^a+\ten_{KLI}^a\ten_{JMI}^a=g_{JK}g_{LM}+g_{LJ}g_{KM}+g_{KL}g_{JM}.$
\end{itemize}
 They reduce
  the ${\bf GL}(\bbR^{n_k})$ group, via $\og(n_k)$, to its subgroup
  $H_k$, where $H_k$ is the irreducible $\sog(3)$ in $\sog(5)$ if
  $k=1$, the irreducible $\sug(3)$ in $\sog(8)$ if $k=2$, the
  irreducible $\spg(3)$ in $\sog(14)$ if $k=4$ and the irreducible
  ${\bf F}_4$ in $\sog(26)$ if $k=8$.

If $k=1,2$ the tensors $\ten_{IJK}^1$ and $\ten_{IJK}^2$ coincide.
If $k=3,4$ they belong to the same $\og(n_k)$ orbit and are related by
  the element ${\rm diag}(1,1,1,1,1,-1,-1,,...,-1)$ of $\og(n_k)$. For
    $k=3,4$ tensors $\ten_{IJK}^1$ and $\ten_{IJK}^2$ are not
    equivalent under the $\sog(n_k)$ action.\label{pb}
\end{proposition}

The above proposition gives examples of a tensor with all the properties
of tensor $\ten$ in dimensions $n=5,8,14$ and $26$. It is remarkable
that it can be proven that these examples exhaust all the
possibilities!

To discuss this statement we need to invoke Elie Cartan's results on
`isoparametric hypersurfaces in spheres' \cite{Cartan}.

\section{Isoparametric hypersurfaces in spheres}

We recall (see. e.g. \cite{Munzner}) that a hypersurface in a real
Riemannian manifold $M$ of constant curvature
is called {\it isoparametric} iff it has constant principal
curvatures. Tuglio Levi-Civita \cite{LC} knew that the number of
such distinct curvatures was at most {\it two} for the Euclidean space
$M=\bbR^3$. The case of $M=\bbR^n$ with arbitrary $n>3$ is similar. It was
shown by Beniamino Segre \cite{Segre} that irrespectively of $n$ the
number of distinct principal curvatures of an isoparametric hypersurface
in $M=\bbR^n$ is at most two. Elie Cartan
\cite{Cartan1} extended this result to the isoparametric
hypersurfaces in the hyperbolic spaces $H^n$ again
showing that in such case the number of possible distinct
principal curvatures is at most two. The situation is quite different
for isoparametric hypersurfaces in spheres $S^n$. In particular Cartan
in Ref. \cite{Cartan} found examples of isoparametric hypersurfaces in
spheres with {\it three} different principal curvatures, each of which
had the same multiplicity. He also introduced
a homogeneous polynomial $$F:\bbR^n\to\bbR$$
of degree $p$ satisfying the differential equations
\begin{eqnarray}
&&{\rm Cii)}\quad\quad\quad  \triangle
  F=0\\\label{carmun1}
&&{\rm Ciii)}\quad\quad\quad |\vec{\nabla} F|^2=p^2~
  \big(~(x^1)^2+
(x^2)^2+...+(x^n)^2~\big)^{(p-1)}.\label{carmun2}
\end{eqnarray}
and proved that all the isoparametric hypersurfaces in $S^{n-1}$ which
have $p$ different constant principal curvatures of the same 
{\it multiplicity} are given
by
\be
S_c=\{x^I\in \bbR^n~|~ F=c~~{\rm and}~~(x^1)^2+
(x^2)^2+...+(x^n)^2=1\},\label{carmun3}
\ee
i.e. that they are the level surfaces of such polynomials $F$
restricted to the sphere.
% Almost 40 years later, ...
%M\"unzner \cite{Munzner} introduced
%a homogeneous polynomial $$F:\bbR^n\to\bbR$$
%of degree $p$ satisfying the differential equations
%\begin{eqnarray*}
%&&{\rm Mii)}\quad\quad\quad  \triangle
%  F=\frac{m_2-m_1}{2~}p^2~\big(~(x^1)^2+
%(x^2)^2+...+(x^n)^2~\big)^{\frac{(p-2)}{2}}\\
%&&{\rm Miii)}\quad\quad\quad |\nabla F|^2=p^2~
%  \big(~(x^1)^2+
%(x^2)^2+...+(x^n)^2~\big)^{(p-1)}.
%\end{eqnarray*}
%In these equations $\nabla$ and $\triangle$ denote, respectively, the
%standard gradient and the standard Laplacian in $\bbR^n$ and $m_2$,
%$m_1$ are real constants. M\"unzner's major result is that
%{\it any} isoparametric hypersurface $S_c$ in $S^{n-1}$ which has $p$
%distinct principal curvatures is obtained from polynomials $F$ by
%Here $c$ is a real number. If all $p$ principal curvatures of the
%hypersurface have the same {\it multiplicities} then $m_1=m_2$ and
%their corresponding polynomial $F$ is {\it harmonic}. This special
%case of M\"unzner's result was already known to Cartan \cite{Cartan},
%who knew that harmonic
%homogeneous polynomials of degree $p$ satisfying condition Miii)
%generate all .

Cartan found {\it all} the homogeneous harmonic polynomials of degree $p=3$
satisfying condition Ciii). On doing that he {\it proved}
\cite{Cartan} that
such polynomials can exist only if $n=n_k=5,8,14,26$. In these four
dimensions he found that the most general form of the polynomials is
$$F=\sum_{I,J,K=1}^{n_k}\ten_{IJK}^ax^Ix^Jx^K,$$
were $\ten_{IJK}^a$ is one of the two tensors appearing in
Proposition \ref{pb}. Writing a generic homogeneous polynomial of
degree $p$ as $F=\ten_{I_1I_2...I_p}x^{I_1}x^{I_2}...x^{I_p}$ we see
that it satisfies Cartan's conditions Cii)-Ciii) iff the {\it
totally symmetric} tensor $\ten_{I_1I_2...I_p}$ satisfies
\begin{eqnarray*}
&&{\rm Cii')}\quad\quad\quad \ten_{JJJ_3J_4...J_p}=0  \\
&&{\rm Ciii')}\quad\quad\quad \ten_{J(J_2J_3...J_p}\ten_{K_2K_3...K_p)J}=g_{(J_2K_2}g_{J_3K_3}...g_{J_pK_p)},
\end{eqnarray*}
where $g_{IJ}={\rm diag}(1,1,...,1)$. Note that in case $p=3$ the
above tensor reduces to $\ten_{IJK}$ and conditions Cii') and Ciii')
become exactly the respective conditions ii) and iii) of Proposition
\ref{pb}. Since $\ten_{IJK}$ is totally symmetric also the condition
i) is satisfied. Thus, Cartan's finding of {\it 
all} isoparametric hypersurfaces with three constant distinct principal
curvatures of the same multiplicity solves our problem of dimensions
in which the tensor $\ten$ may exist. Summarizing we have the
following theorem, which is a reformulation of the above mentioned
Cartan's results.
\begin{theorem}
An $\bbR^n$ with the standard Euclidean metric $g_{IJ}$ admits
a tensor $\ten_{IJK}$ with the properties
\begin{itemize}
\item[i)] %it was totally symmetric:
$\ten_{IJK}=\ten_{(IJK)}$,
\vskip1mm
\item[ii)]%it was trace free:
$\ten_{IJJ}=0$,
\vskip1mm
\item[iii)] $\ten_{JKI}\ten_{LMI}+\ten_{LJI}\ten_{KMI}+\ten_{KLI}\ten_{JMI}=g_{JK}g_{LM}+g_{LJ}g_{KM}+g_{KL}g_{JM}$
%it satisfied the following condition:
\end{itemize}
if and only if $n=n_k=3k+2$ for $k=1,2,4,8$. Modulo the action of the
$\sog(n_k)$ group all such tensors are given by Proposition \ref{pb}.
\end{theorem}

\section{Representations of $\sug(3)$, $\spg(3)$ and
  ${\bf F}_4$.}

It is known that there are {\it real} irreducible representations of the
group $\sug(3)$ in dimensions:
$$1,8,20,27,70.$$
Also, there are {\it real} irreducible representations of the
group $\spg(3)$ in dimensions:
$$1,14,21,70,84,90,126,189,512,525$$
and there are {\it real} irreducible representations of the
group ${\bf F}_4$ in dimensions:
$$1,26,52,273,324,1053,1274,4096,8424.$$
To see how these representations appear we consider a vector space
$\bbR^{n_k}$, $n_k=5,8,14,26$ equipped with the Riemannian metric
$g$ and the corresponding tensor $\ten^1_{IJK}$ of Proposition
\ref{pb}. As we know the stabilizer $H_k$ of $\ten^1$ is a subgroup
of $\sog(n_k)$, which when $n_k$ is 5, 8, 14 and 26 is,
respectively, $\sog(3),\sug(3),\spg(3)$ and ${\bf F}_4$. Now, the
tensor $\ten^1$ can be used to decompose the tensor product
representation $\bgt^2\bbR^{n_k}$ of the group $H_k$ onto the {\it
real} irreducible components as follows. First, we define an
endomorphism \be\hat{\ten}:{\textstyle \bigotimes}^2\bbR^{n_k}
\longrightarrow{\textstyle\bigotimes}^2\bbR^{n_k},\label{hatten}\ee
$$W^{IK}\stackrel{\hat{\ten}}{\longmapsto} 4~\ten^1_{IJM}\ten^1_{KLM}
    W^{JL},$$
which preserves the decomposition
    $\bigotimes^2\bbR^{n_k}=\bgw^2\bbR^{n_k}\oplus\bgs^2\bbR^{n_k}$.
Second, we look for its eigenspaces, which surely are $H_k$ invariant.
    We have the following proposition.
\begin{proposition} \label{pr:so3rep}~\\
1) If $n_k=5$ then
$$\bgt^2\bbR^{n_k} = {^5\bgw^2_3} \oplus {^5\bgw^2_7} \oplus {^5\bgs^2_1} \oplus {^5\bgs^2_5} \oplus
  {^5\bgs^2_9},$$
where
\begin{eqnarray*}
  &&^5\bgs^2_1       = \{~S\in\bgt^2\bbR^5~|~ \hat{\ten} (S)= 14 \cdot S~\} = \{S=\lambda
  \cdot g,\;\lambda\in\bbR~\}, \\
 && ^5\bgw^2_3 = \{ ~F\in\bgt^2 \bbR^5~|~  \hat{\ten} (F) = 7 \cdot F~\}~
=~\soa(3),\\
 && ^5\bgs^2_5 = \{~ S\in\bgt^2\bbR^5 ~| ~ \hat{\ten}(S) = -3 \cdot S~\}, \\
 && ^5\bgw^2_7 = \{ ~F\in\bgt^2\bbR^5 ~| ~  \hat{\ten}(F) = -8 \cdot F~\},\\
 &&^5\bgs^2_9 = \{~ S\in\bgt^2\bbR^5 ~| ~ \hat{\ten}(S) = 4\cdot S~\}.
\end{eqnarray*}
The \emph{real} vector spaces $^5\bgw^2_i\subset\bgw^2\bbR^5$ and
$^5\bgs^2_j\subset\bgs^2\bbR^5$ of respective dimensions $i$ and $j$ are irreducible representations of the
group $\sog(3)$.\\
2) If $n_k=8$ then
$$\bgt^2\bbR^{n_k} = {^8\bgw^2_8} \oplus {^8\bgw^2_{20}} \oplus {^8\bgs^2_1} \oplus {^8\bgs^2_8} \oplus
  {^8\bgs^2_{27}},$$
where
\begin{eqnarray*}
  &&^8\bgs^2_1       = \{~S\in\bgt^2\bbR^8~|~ \hat{\ten} (S)= 20 \cdot S~\} = \{S=\lambda
  \cdot g,\;\lambda\in\bbR~\}, \\
 && ^8\bgw^2_8 = \{ ~F\in\bgt^2 \bbR^8~|~  \hat{\ten} (F) = 10 \cdot F~\}~
=~\sua(3),\\
 && ^8\bgs^2_8 = \{~ S\in\bgt^2\bbR^8 ~| ~ \hat{\ten}(S) = -6 \cdot S~\}, \\
 && ^8\bgw^2_{20} = \{ ~F\in\bgt^2\bbR^8 ~| ~  \hat{\ten}(F) = -8 \cdot F~\},\\
 &&^8\bgs^2_{27} = \{~ S\in\bgt^2\bbR^8 ~| ~ \hat{\ten}(S) = 4\cdot S~\}.
\end{eqnarray*}
The \emph{real} vector spaces $^8\bgw^2_i\subset\bgw^2\bbR^8$ and
$^8\bgs^2_j\subset\bgs^2\bbR^8$ of respective dimensions $i$ and $j$ are irreducible representations of the
group $\sug(3)$. The representations $\bgw^2_8$
and $\bgs^2_8$ are \emph{equivalent}.\\
3) If $n_k=14$ then
$$\bgt^2\bbR^{n_k} = {^{14}\bgw^2_{21}} \oplus {^{14}\bgw^2_{70}} \oplus {^{14}\bgs^2_1} \oplus {^{14}\bgs^2_{14}} \oplus
  {^{14}\bgs^2_{90}},$$
where
\begin{eqnarray*}
  &&^{14}\bgs^2_1       = \{~S\in\bgt^2\bbR^{14}~|~ \hat{\ten} (S)= 32 \cdot S~\} = \{S=\lambda
  \cdot g,\;\lambda\in\bbR~\}, \\
 && ^{14}\bgw^2_{21} = \{ ~F\in\bgt^2 \bbR^{14}~|~  \hat{\ten} (F) = 16 \cdot F~\}~
=~\spa(3),\\
 && ^{14}\bgs^2_{14} = \{~ S\in\bgt^2\bbR^{14} ~| ~ \hat{\ten}(S) = -12 \cdot S~\}, \\
 && ^{14}\bgw^2_{70} = \{ ~F\in\bgt^2\bbR^{14} ~| ~  \hat{\ten}(F) = -8 \cdot F~\},\\
 &&^{14}\bgs^2_{90} = \{~ S\in\bgt^2\bbR^{14} ~| ~ \hat{\ten}(S) = 4\cdot S~\}.
\end{eqnarray*}
The \emph{real} vector spaces $^{14}\bgw^2_i\subset\bgw^2\bbR^{14}$ and
$^{14}\bgs^2_j\subset\bgs^2\bbR^{14}$  of respective dimensions $i$ and $j$ are irreducible representations of the
group $\spg(3)$.\\
4) If $n_k=26$ then
$$\bgt^2\bbR^{n_k} = {^{26}\bgw^2_{52}} \oplus {^{26}\bgw^2_{273}} \oplus {^{26}\bgs^2_1} \oplus {^{26}\bgs^2_{26}} \oplus
  {^{26}\bgs^2_{324}},$$
where
\begin{eqnarray*}
  &&^{26}\bgs^2_1       = \{~S\in\bgt^2\bbR^{26}~|~ \hat{\ten} (S)= 56 \cdot S~\} = \{S=\lambda
  \cdot g,\;\lambda\in\bbR~\}, \\
 && ^{26}\bgw^2_{52} = \{ ~F\in\bgt^2 \bbR^{26}~|~  \hat{\ten} (F) = 28 \cdot F~\}~
=~{\frak f}_4,\\
 && ^{26}\bgs^2_{26} = \{~ S\in\bgt^2\bbR^{26} ~| ~ \hat{\ten}(S) = -24 \cdot S~\}, \\
 && ^{26}\bgw^2_{273} = \{ ~F\in\bgt^2\bbR^{26} ~| ~  \hat{\ten}(F) = -8 \cdot F~\},\\
 &&^{26}\bgs^2_{324} = \{~ S\in\bgt^2\bbR^{26} ~| ~ \hat{\ten}(S) = 4\cdot S~\}.
\end{eqnarray*}
The \emph{real} vector spaces $^{26}\bgw^2_i\subset\bgw^2\bbR^{26}$ and
$^{26}\bgs^2_j\subset\bgs^2\bbR^{26}$  of respective dimensions $i$ and $j$
are irreducible representations of the
group ${\bf F}_4$.
\end{proposition}
\begin{remark}
According to this proposition we may identify spaces $\bbR^{n_k}$ with
the representation spaces $^{n_k}\bgs^2_{n_k}$ corresponding to the
eigenvalues $2-n_k$ of $\hat{\ten}$.
Noting that the dimension of the group $H_k$ is $$\dim
H_k=4k-1+(k-1)\log_2 k,\quad\quad\quad k=1,2,4,8$$
and introducing $s_k=\tfrac{9}{2}k(k+1)$, we see that
\begin{itemize}
\item the eigenvalues of $\hat{\ten}$ corresponding to spaces
$^{n_k}\bgs^2_1$ are $4+2n_k$,
\item the eigenvalues corresponding to
spaces $^{n_k}\bgw^2_{\dim H_k}$ are $2+n_k$,
\item the eigenvalues corresponding to
spaces $^{n_k}\bgw^2_{(s_k+1-\dim H_k)}$ are
\emph{always} $-8$,
\item the eigenvalues corresponding to spaces
$^{n_k}\bgs^2_{s_k}$ are \emph{always} $+4$.
\end{itemize}
We also note that we may
identify the Lie algebras $\frak{h}_k$ of $H_k$ with the representations $^{n_k}\bgw^2_{\dim H_k}$.
\end{remark}

\section{$H_k$ structures on Riemannian manifolds}\label{hkst}
\begin{definition}
An $H_k$ structure on a $n_k$-dimensional Riemannian manifold $(M,g)$
is a structure defined by means of a rank 3 tensor field $\ten$
satisfying
\begin{enumerate}
\item[i)] $\ten_{IJK}=\ten_{(IJK)}$,
\item[ii)]%it was trace free:
$\ten_{IJJ}=0$,
\vskip1mm
\item[iii)]
  $\ten_{JKI}\ten_{LMI}+\ten_{LJI}\ten_{KMI}+\ten_{KLI}\ten_{JMI}=g_{JK}g_{LM}+g_{LJ}g_{KM}+g_{KL}g_{JM}.$
\end{enumerate}
\end{definition}

\begin{definition}
Two $H_k$ structures $(M,g,\ten)$ and $(\bar{M},\bar{g},\bar{\ten})$
defined on
two respective $n_k$-manifolds $M$ and $\bar{M}$ are (locally) \emph{equivalent}
iff there exists a (local) diffeomorphism $\phi:M\to \bar{M}$ such
that
$$\phi^*(\bar{g})=g\quad\quad{\rm and}\quad\quad \phi^*(\bar{\ten})=\ten.$$
If $\bar{M}=M$, $\bar{g}=g$, $\bar{\ten}=\ten$ the equivalence $\phi$ is
called a (local) \emph{symmetry} of $(M,g,\ten)$. The group of (local)
symmetries is called a \emph{symmetry group of} $(M,g,\ten)$.
\end{definition}

\noindent As we now the tensor field $\ten$ reduces the structure
group of the bundle of orthonormal frames over $M$ to one of the
groups $H_k$ of Proposition \ref{pb}. We also know that the Lie
algebra $\mathfrak{h}_k$ of $H_k$ is isomorphic to
$\mathfrak{h}_k\simeq{^{n_k}\Lambda}^2_{\dim
H_k}\subset\otimes^2\bbR^{n_k}$ of Proposition \ref{pr:so3rep}.
Thus, at each point, every element $F$ of the Lie algebra
$\mathfrak{h}_k$ may be considered to be an endomorphism of
$\bbR^{n_k}$. This defines an element $$f=\exp(F)\in H_k\subset
\sog(n_k)\subset{\bf GL}(n_k,\bbR)$$ and, point by point, induces
the natural action $\rho(f)$ of the group $H_k$ on the vector-valued
1-forms
$$\theta=(\theta^1,\theta^2,\theta^3,\dots, \theta^{n_k} )\in\bbR^{n_k}\otimes\Om^1(M)$$ by:
\be \theta\mapsto \thet=\rho(f)(\theta)=f\cdot\theta.\label{lifted}
\ee This, enables for local description of a $H_k$ structure on $M$
by means of a coframe \be
\theta=(\theta^I)=(\theta^1,\theta^2,\theta^3,\dots, \theta^{n_k}
)\label{coframe} \ee on $M$, given up to the $H_k$ transformations
(\ref{lifted}).

For such a class of coframes the Riemannian metric $g$ is
$$
g=\theta_1^2+\theta_2^2+\theta_3^2+\dots+\theta_{n_k}^2,
$$
and the tensor $\ten$, reducing the structure group from $\sog(n_k)$ to
$H_k$, is
\begin{equation}
\label{tform}
\ten=\ten^1_{IJK}\theta^I\theta^J\theta^K,
\end{equation}
where $\ten^1$ is defined in Proposition \ref{pb}.
\begin{definition}
\label{df:adapted}
An orthonormal coframe
$(\theta^1,\theta^2,\theta^3,\dots,\theta^{n_k})$ in which
the tensor $\ten$ of an $H_k$ structure $(M,g,\ten)$ is of the form
(\ref{tform}) is called a \emph{coframe adapted to $(M,g,\ten)$}, an
\emph{adapted coframe}, for short.
\end{definition}

 Given an $H_k$
structure as above, we consider an {\it arbitrary}
  $\mathfrak{h}_k$-valued connection on $M$. This may be locally
represented by means of an $\mathfrak{h}_k$-valued 1-form $\Gamma$
given by \be \Gamma=(\Gamma^I_{~J})=\gamma^\alpha
E_\alpha,\quad\quad\alpha=1,2,\dots,\dim H_k\label{charcon} \ee
where $\gamma^\alpha$ are 1-forms on $M$ and for each $\alpha$ the
  symbols
  $E_\alpha=({E_\alpha}^I_{~J})$ denote constant $n_k\times n_k$-matrices
  which form a basis of the Lie algebra $\mathfrak{h}_k$. The explicit
  expressions for $E_\alpha$ are presented in the Appendix.
The connection $\Gamma$, having values in
$\mathfrak{h}_k\subset\soa(n_k)$, is
  necessarily metric. Via the Cartan structure equations,
\be \der\theta^I+\Gamma^I_{~J}\dz\theta^J=T^I
%=\frac{1}{2}T^i_{~jk}\theta^j\dz\theta^k,
\label{chartor} \ee \be
\der\Gamma^I_{~J}+\Gamma^I_{~K}\dz\Gamma^K_{~J}=R^I_{~J},
%=\frac{1}{2}K^i_{~jkm}\theta^k\dz\theta^m,
\label{charcurv} \ee it defines the torsion 2-form $T^I$ and the
$\mathfrak{h}_k$-curvature 2-form $R^I_{~J}$. Using these forms we
define the torsion tensor
$T^I_{~JK}\in(\bbR^{n_k}\otimes\bigwedge^2\bbR^{n_k})$ and the
$\mathfrak{h}_k$-curvature tensor
$r^\alpha_{~JK}\in(\mathfrak{h}_k\otimes\bigwedge^2\bbR^{n_k})$,
respectively, by
$$T^I= \frac{1}{2}T^I_{~JK}\theta^J\dz\theta^K$$
and \be R^I_{~J}=\frac12 r^\alpha_{KL}\theta^K\dz\theta^L
{E_\alpha}^I_{~J}.\label{ccu} \ee The connection satisfies the first
Bianchi identity \be R^I_{~J}\dz\theta^J=DT^I\label{bianchi1} \ee
and the second Bianchi identity \be DR^I_{~J}=0,\label{bianchi2} \ee
with the covariant differential defined by
$$DT^I=\der T^I+\Gamma^I_{~J}\dz T^J,\quad\quad\quad DR^I_{~J}=\der
R^I_{~J}+\Gamma^I_{~K}\dz R^K_{~J}-R^I_{~K}\dz\Gamma^K_{~J}.$$

Since the $H_k$ preserves $g$ and $\ten$ we have the following
proposition.
\begin{proposition}\label{oco}
Every $\mathfrak{h}_k$-valued connection $\Gamma$ of (\ref{charcon})
is metric
$$\stackrel{\Gamma}{\nabla}_v(g)\equiv 0$$
and preserves tensor $\ten$
$$\stackrel{\Gamma}{\nabla}_v(\ten)\equiv 0\quad\quad\quad
\forall v\in{\rm T}M.$$
\end{proposition}
\noindent

\section{Characteristic connection}\label{cco}
In this section we consider $H_k$ structures $(M,g,\ten)$ with
Levi-Civita connection $\lc\in\soa(n_k)\otimes\bbR^{n_k}$ {\it
uniquely} decomposable according to \be \lc=\Gamma+\tfrac{1}{2}T,
\label{rozk} \ee where $\Gamma\in \mathfrak{h}_k\otimes\bbR^{n_k}$
and $T\in \bgw^3\bbR^{n_k}$.

Such $H_k$ structures are interesting, since for them, contrary to
the generic case, the decomposition (\ref{rozk}) defines a unique
$\mathfrak{h}_k$-valued connection $\Gamma$. Moreover, given the
unique decomposition (\ref{rozk}), we may rewrite the Cartan
structure equations
$$\der \theta^I+\lc~^I_{~J}\wedge\theta^J=0$$
for the Levi-Civita connection $\lc$ into the form
$$\der
\theta^I+\Gamma^I_{~J}\wedge\theta^J=\tfrac{1}{2}T^I_{~JK}\theta^J\wedge\theta^K$$
and to interpret $T$ as the {\it totally skew symmetric
  torsion} of $\Gamma$.
\begin{definition}
A $\mathfrak{h}_k$-valued connection $\Gamma$ of a $H_k$ structure
$(M,g,\ten)$ admitting the unique decomposition
$$\lc=\Gamma+\tfrac{1}{2}T,\quad{\rm with}\quad\Gamma\in \mathfrak{h}_k\otimes\bbR^{n_k}\quad{\rm and}
\quad T\in \bgw^3\bbR^{n_k}$$ is called the characteristic
connection.
\end{definition}
Since $\Gamma\in\mathfrak{h}_k\otimes\bbR^{n_k}$ and
$T\in\bgw^3\bbR^{n_K}$ it is obvious from (\ref{rozk}) that the
Levi-Civita connection $\lc$ of $H_k$ structures which admit
characteristic connections must satisfy
\be\lc\in[\mathfrak{h}_k\otimes\bbR^{n_k}]+\bgw^3\bbR^{n_k}.\label{kole}\ee
Moreover, since
$\dim(\mathfrak{h}_k\otimes\bbR^{n_k})+\dim(\bgw^3\bbR^{n_k})
<\dim(\soa(n_k)\otimes\bbR^{n_k})$ then it is obvious that the
unique decomposition (\ref{rozk}) is not possible for all $H_k$
structures. Our aim now is to characterize $H_k$ structures
admitting characteristic connection.

Following \cite{bobi} we introduce the following definition.

\begin{definition}
An $H_k$ structure $(M,g,\ten)$ is called \emph{nearly integrable}
iff \be (\stackrel{\scriptscriptstyle{LC}}{\nabla}_v
\ten)(v,v,v)\equiv 0\label{ni} \ee for the Levi-Civita connection
$\stackrel{\scriptscriptstyle{LC}}{\nabla}$ and for every vector field $v$ on $M$..
\end{definition}
The condition (\ref{ni}), when written in an adapted coframe
(\ref{coframe}), is \be \lc_{M(JI}~\ten_{KL)M}\equiv
0,\label{nearly} \ee where $\lc_{MJ}=\lc_{MJK}\theta^K$ denotes the
$\soa(n_k)$-valued 1-form corresponding to the Levi-Civita
connection $\stackrel{\scriptscriptstyle{LC}}{\nabla}$. This
motivates an introduction of the map
$$\ten':\bgw^2\bbR^{n_k}\otimes\bbR^{n_k}\mapsto \bgs^4\bbR^{n_k}$$
such that \beq
\ten'(\lc)_{IJKL}&=&12\lc_{M(JI}~\ten_{KL)M}\nonumber\\
&=&\lc_{MJI}~\ten_{MKL}+\lc_{MKI}~\ten_{JML}+\lc_{MLI}~\ten_{JKM}\nonumber\\
&+&\lc_{MIJ}~\ten_{MKL}+\lc_{MKJ}~\ten_{IML}+\lc_{MLJ}~\ten_{IKM}\label{tprim}\\
&+&\lc_{MIK}~\ten_{MJL}+\lc_{MJK}~\ten_{IML}+\lc_{MLK}~\ten_{IJM}\nonumber\\
&+&\lc_{MIL}~\ten_{MJK}+\lc_{MJL}~\ten_{IMK}+\lc_{MKL}~\ten_{IJM}.\nonumber
\eeq Comparing this with (\ref{nearly}) we have the following
proposition.
\begin{proposition}
An $H_k$ structure $(M,g,\ten)$ is nearly integrable if and only if
its Levi-Civita connection $\lc\in\ker \ten'$.
\end{proposition}
It is worthwhile to note that each of the last four rows of
(\ref{tprim}) resembles the l.h.s. of equality
$$
X_{MJ}\ten_{MKL}+X_{MK}\ten_{JML}+X_{ML}\ten_{JKM}=0
$$
satisfied by every matrix $X\in\mathfrak{h}_k={^{n_k}}\bgw^2_{\dim
H_k}$. Thus, $\mathfrak{h}_k\otimes\bbR^{n_k}\subset\ker \ten'.$ Due
to the first equality in (\ref{tprim}) we also have
$\bgw^3\bbR^{n_k}\subset\ker \ten'.$ This proves the following
Lemma.
\begin{lemma}
Since
$$\mathfrak{h}_k\otimes\bbR^{n_k}\subset\ker \ten'\quad\quad{\rm
and}\quad\quad\bgw^3\bbR^{n_k}\subset\ker \ten'$$ then
$$([\mathfrak{h}_k\otimes\bbR^{n_k}]+\bgw^3\bbR^{n_k})\subset\ker
\ten'.$$
\end{lemma}
Thus, comparing this with (\ref{kole}) we have the following
Proposition.
\begin{proposition}
Among all $H_k$ structures only the nearly integrable ones
\emph{may} admit characteristic connection.
\end{proposition}

It is known \cite{bobi} that if $n_k=5$ then the nearly
integrability condition is also sufficient for the existence of
a characteristic connection. To see that it is
no longer true for all $n_k$ we need to see how the intersections
$[\mathfrak{h}_k\otimes\bbR^{n_k}]\cap\bgw^3\bbR^{n_k}$ and the
algebraic sums $[\mathfrak{h}_k\otimes\bbR^{n_k}]+\bgw^3\bbR^{n_k}$
depend on the dimension $n_k$. After some algebra we arrive at the
following proposition.
\begin{proposition}~\label{chp}\\
\vspace{-.5cm}
\begin{itemize}
\item If $n_k=5$ or $n_k=14$ then
$$\ker\ten'=[\mathfrak{h}_k\otimes\bbR^{n_k}]\oplus\bgw^3\bbR^{n_k}.$$
\item If $n_k=8$ then
  $$\ker\ten'=[\sua(3)\otimes\bbR^{8}]+\bgw^3\bbR^{8}\quad and\quad[\sua(3)\otimes\bbR^{8}]\cap\bgw^3\bbR^{8}={^8\bgs}^2_1.$$
\item If $n_k=26$ then
  $$\ker\ten'=[\mathfrak{f}_4\otimes\bbR^{26}]\oplus\bgw^3\bbR^{26}\oplus{^{26}\bgw}^2_{52}.$$
\item In particular, for $n_k=5,14$ and $26$ we have $[\mathfrak{h}_k\otimes\bbR^{n_k}]\cap\bgw^3\bbR^{n_k}=\{0\}$.
\end{itemize}
\end{proposition}
This implies the following Theorem.
\begin{theorem}
In dimensions $n_k=5$ and $n_k=14$ the necessary and sufficient
condition for a $H_k$ structure $(M,g,\ten)$ to admit a
characteristic connection is that $(M,g,\ten)$ is nearly integrable
$$(\stackrel{\scriptscriptstyle{LC}}{\nabla}_v
\ten)(v,v,v)\equiv 0.$$
\end{theorem}
Proposition \ref{chp} also implies that the nearly integrable $H_k$
structures in dimension $n_k=8$ admit decomposition (\ref{rozk}).
However, in this dimension condition (\ref{rozk}) determins the 
connection $\Gamma$ and the torsion
$T$ up to an additional freedom. Due to the 1-dimensional intersection
$[\sua(3)\otimes\bbR^{8}]\cap\bgw^3\bbR^{8}={^8\bgs}^2_1$ we see
that in such case there is a 1-parameter family of connections
$\Gamma(\lambda)\in\sua(3)\otimes\bbR^{8}$ and 1-parameter family of
skew symmetric torsions $T(\lambda)\in\bgw^3\bbR^{8}$ such that \be
\lc=\Gamma(\lambda)+\tfrac12 T(\lambda).\label{odla} \ee It is clear
that for $n_k=8$, the requirement (\ref{rozk}) uniquely determines
$\Gamma\in \sua(3)\otimes\bbR^{8}$ and $T\in \bgw^3\bbR^{8}$ only if
we restrict ourselves to the nearly integrable $\sug(3)$ structures
for which the Levi-Civita connection $\lc$ is in the 118-dimensional
space ${^8\mathcal{V}}$ such that
${^8\mathcal{V}}\oplus{^8\bgs}^2_1=\ker\ten'$. It follows that this
space has the following decomposition under the $\sug(3)$ action
${^8\mathcal{V}}=2~{^8\bgs^2_{27}} \oplus 2~{^8\bgw^2_{20}} \oplus
3~{^8\bgs^2_8}$. It is convenient to extend this notation and to
introduce vector spaces ${^{n_k}\mathcal{V}}$ to be
  subspaces of $\ker\ten'$ such that:
\begin{eqnarray*}
{^{n_k}\mathcal{V}}&=&\ker\ten'\quad{\rm for}\quad n_k=5,14\\
{^8\mathcal{V}}&=&2~{^8\bgs^2_{27}} \oplus 2~{^8\bgw^2_{20}} \oplus
3~{^8\bgs^2_8}\varsubsetneq\ker\ten',\\
{^{26}\mathcal{V}}&=&[\mathfrak{f}_4\otimes\bbR^{26}]\oplus\bgw^3\bbR^{26}\varsubsetneq\ker\ten'.
\end{eqnarray*}

Using these we have the following definition

\begin{definition}
An $H_k$ structure $(M,g,\ten)$ is called \emph{restricted nearly
integrable} iff its Levi-Civita connection
$\lc\in{^{n_k}\mathcal{V}}$.
\end{definition}

\begin{remark}
Note that for $n_k=5$ or $n_k=14$ the term: {\it restricted nearly
integrable} is the same as: {\it nearly integrable}.
\end{remark}
Looking again at Proposition \ref{chp} we see that the above
restriction for the nearly integrable $\sug(3)$ or ${\bf F}_4$
structures in respective dimensions $n_k=8$ and $n_k=26$ is
precisely the one that gives the sufficient conditions for the
existence and uniqueness of the characteristic connection.
Summarizing we have the following theorem.
\begin{theorem}
A necessary and sufficient condition for a $H_k$ structure
$(M,g,\ten)$ to admit a characteristic connection is that this
structure is restricted nearly integrable.
\end{theorem}
\begin{remark}
Note, that if $n_k=5$ then, out of the {\it a priori} 50 independent
components of the Levi-Civita connection $\lc$, the (restricted)
nearly integrable condition (\ref{rozk}) excludes 25. Thus,
heuristically, the (restricted) nearly integrable $\sog(3)$
structures constitute `a half' of all the possible $\sog(3)$
structures in dimension 5.

If $n_k=8$ the Levi-Civita connection has 224 components. The
restricted nearly integrable condition reduces it to 118. For
$n_k=14$ these numbers reduce from 1274 to 658. For $n_k=26$ the
reduction is from 8450 to 3952.
\end{remark}

\section{Classification of the restricted nearly integrable $H_k$
structures} We classify the possible types of the restricted nearly
integrable $H_k$ structures according to the $H_k$ irreducible
decompositions of the spaces $\bgw^3\bbR^{n_k}$ in which the torsion
$T$ of their characteristic connections live. Using a computer
algebra package for Lie group computations `LiE' \cite{lie} we
easily arrive at the following proposition.
\begin{proposition}\label{clastor}
Let $(M,g,\ten)$ be a restricted nearly integrable $H_k$ structure.
The $H_k$ irreducible decomposition of the skew symmetric torsion
$T$ of the characteristic connection for $(M,g,\ten)$ is given by:
\begin{itemize}
\item
$T\in {^5\bgw^2_7}\oplus{^5\bgw^2_3},$ \hspace{3.5cm}for $n_k=5$,
\item
$T\in
  {^8\bgs^2_{27}}\oplus{^8\bgw^2_{20}}\oplus{^8\bgs^2_{8}}\oplus{^8\bgs^2_{1}},$
  \hspace{.87cm}for $n_k=8$,
\item $T\in
  {^{14}V_{189}}\oplus{^{14}V_{84}}\oplus{^{14}\bgw^2_{70}}\oplus{^{14}\bgw^2_{21}},$
 \hspace{.36cm}for $n_k=14$,
\item
$T\in
{^{26}V_{1274}}\oplus{^{26}V_{1053}}\oplus{^{26}\bgw^2_{273}},$
 \hspace{1.1cm}for $n_k=26$.
\end{itemize}
Here ${^{n_k}V_j}$ denotes irreducible $j$-dimensional representations of 
$H_k$ which were not present in the $H_k$ decomposition of $\bigotimes^2\bbR^{n_k}$.
\end{proposition}

This provides an analog of the Gray-Hervella \cite{gray}
classification for the restricted nearly integrable $H_k$
structures.

We close this section with a remark on possible types of the
curvature $R$ of the characteristic connections.

\begin{remark}
In the below formulae $^{n_k}V_j$ denote $j$-dimensional irreducible
representation space for $H_k$ which did not appear in Proposition
\ref{pr:so3rep}.
\begin{itemize}
\item If $n_k=5$ then $R\in{^5\bgs^2_9} \oplus{^5\bgw^2_7}
\oplus   2~ {^5\bgs^2_5} \oplus{^5\bgw^2_3} \oplus {^5\bgs^2_1}$
\item If $n_k=8$ then $R\in{^8 V_{70}}\oplus 3~{^8\bgs^2_{27}} \oplus
2~{^8\bgw^2_{20}} \oplus 4~{^8\bgs^2_8} \oplus {^8\bgs^2_1}$.
\item If $n_k=14$ then
$R\in{^{14}V_{525}}\oplus{^{14}V_{512}}\oplus2~{^{14}V_{189}}\oplus{^{14}V_{126}}\oplus 2~{^{14}\bgs^2_{90}} \oplus2~{^{14}\bgw^2_{70}}
\oplus{^{14}\bgw^2_{21}} \oplus2~{^{14}\bgs^2_{14}} \oplus
      {^{14}\bgs^2_1}$.
\item If $n_k=26$ then $R\in{^{26}V_{8424}}\oplus{^{26}V_{4096}}\oplus{^{26}V_{1274}}\oplus{^{26}V_{1053}}\oplus{^{26}{V'}_{1053}}\oplus 2~{^{26}\bgs^2_{324}} \oplus{^{26}\bgw^2_{273}}
\oplus{^{26}\bgw^2_{52}} \oplus{^{26}\bgs^2_{26}} \oplus
  {^{26}\bgs^2_1}$.
\end{itemize}
Note that, due to the restricted nearly integrability condition, it
is rather unlikely that $R$ may attain values in all of the above
irreducible parts.
\end{remark}
\section{Dimensions 12, 18, 28, 30, 40, 54, 64 and 112; the `exceptional' 8 and 32}
\subsection{Torsionless models} It is obvious that the simplest restricted nearly
integrable $H_k$ structures have the characteristic connection
$\Gamma$ with vanishing torsion $T\equiv 0$. For them
$$\lc=\Gamma\in\mathfrak{h}_k\otimes\bbR^{n_k},$$
hence their Riemannian holonomy group is reduced from $\sog(n_k)$ to
the group $H_k$. Since $H_k\subset\sog(n_k)$ in the respective
dimensions $n_k=5,8,14$ and $26$ are not present in the Berger list of
the Riemannian holonomy groups \cite{ber}, the only possible
restricted nearly integrable $H_k$ structures with $T\equiv 0$ must
be locally isometric to the symmetric spaces $M=G_k/H_k$. The Lie
group $G_k$ appearing here must to have dimension $\dim G_k=n_k+\dim
H_k$. Looking at Cartan's list \cite{car} of the irreducible
symmetric spaces (see e.g. \cite{bes} pp. 312-317) we have the
following Theorem.
\begin{theorem}\label{torless}
All $H_k$ structures with vanishing torsion are locally isometric to
one of the symmetric spaces
$$M=G_k/H_k,$$ where the possible Lie groups $G$ are given in the following
table:\\
\begin{tabular}{|l|l|l|l|l|}
\hline
 $\dim M$&{\rm Group} $H_k$&{\rm Group} $G_k$&{\rm Group} $G_k$&{\rm Group} $G_k$\\
\hline\hline $n_k=5$&$\sog(3)$&
$\sug(3)$&$\sog(3)\times_\rho\mathbb{R}^5$&${\bf
SL}(3,\bbR)$\\
\hline $n_k=8$&$\sug(3)$&
$\sug(3)\times\sug(3)$&$\sug(3)\times_\rho\mathbb{R}^8$&${\bf
SL}(3,\bbC)$\\
\hline $n_k=14$&$\spg(3)$&
$\sug(6)$&$~\spg(3)\times_\rho\mathbb{R}^{14}$&${\bf
SU}^*(6)\simeq{\bf SL}(3,\bbH)$\\
\hline $n_k=26$&${\bf F}_4$&${\bf E}_6$&$~\,~~~~\quad{\bf
F}_4\times_\rho\mathbb{R}^{26}$&${\bf E}^{-26}_6\simeq{\bf
SL}(3,\bbO)$\\
\hline
\end{tabular}\\
Here $\rho$ is the irreducible representation of $H_k$ in
$\bbR^{n_k}$.
\end{theorem}
\begin{remark}
Let $\mathfrak{g}_k$ be the Lie algebra of the group $G_k$ of
Theorem \ref{torless}. We note that since the torsionless models for
the $H_k$ structures are the symmetric spaces $M=G_k/H_k$, then
arbitrary restricted nearly integrable $H_k$ structures may be
analyzed in terms of a {\it Cartan $\mathfrak{g}_k$-valued
connection} on the Cartan bundle $H_k\to P\to M$. In such a langauge
the {\it torsionless} models with respect to the $\mathfrak{h}_k$ connection are simply the {\it flat} models for the corresponding Cartan
$\mathfrak{g}_k$-valued connection on $P$.
\end{remark}
\begin{remark}
According to \cite{wolf} the manifold $M=\sug(3)/\sog(3)$ is a
unique irreducible Riemannian symmetric space $M=G/H$ with the
property that $({\rm rank} G-{\rm rank} H)=1$ and that $M$ is not
isometric to an odd dimensional real Grassmann manifold. It is
interesting to note \cite{frie} (see \cite{wolf} p. 324) that the
other {\it compact} torsionless $H_k$ structures correspond to
manifolds $M=\sug(3)$, $M=\sug(6)/\spg(3)$ and $M={\bf E}_6/{\bf
F}_4$, which are examples of a very few irreducible symmetric
Riemannian manifolds $M=G/H$ with $({\rm rank} G-{\rm rank} H)=2$.
\end{remark}
\subsection{The `magic square'}\label{ms}
We now concentrate on the Lie algebras $\mathfrak{h}_k$ and $\mathfrak{g}_k$
corresponding to groups $H_k$ and $G_k$ appearing in the second and
the third column of the table included in Theorem \ref{torless}. We note 
that these Lie algebras constitute the first two columns of the `magic
square' \cite{tit,mag}:\\
\begin{center}
\begin{tabular}{|c|c|c|c|}
\hline
$\soa(3)$&$\sua(3)$&$\spa(3)$&$\mathfrak{f}_4$\\
\hline
$\sua(3)$&$2\sua(3)$&$\sua(6)$&$\mathfrak{e}_6$\\
\hline
$\spa(3)$&$\sua(6)$&$\soa(12)$&$\mathfrak{e}_7$\\
\hline
$\mathfrak{f}_4$&$\mathfrak{e}_6$&$\mathfrak{e}_7$&$\mathfrak{e}_8$\\
\hline
\end{tabular}~~.
\end{center}
In agreement with the previous notation let us denote by $H_k$, $G_k$, $\mathcal{G}_k$
and $\tilde{\mathcal{G}}_k$ the compact Lie groups corresponding to the
Lie algebras of the first, the second, the third and the fourth respective columns
of the magic square. Since $G_k/H_k$ are the torsionless compact models for 
$H_k$ geometries, it may seem reasonable to consider spaces 
$\mathcal{G}_k/G_k$ and $\tilde{\mathcal{G}}_k/\mathcal{G}_k$ as the torsionless models for new special Riemannian geometries with a characteristic connection. 
Unfortunately the homogeneous spaces $\mathcal{G}_k/G_k$ and
$\tilde{\mathcal{G}}_k/\mathcal{G}_k$  are {\it reducible}. However, if we replace the second
column in the magic square by 
\begin{center} 
\begin{tabular}{|c|}
\hline
$\sua(3)\oplus\bbR$\\
\hline
$2\sua(3)\oplus\bbR$\\
\hline
$\sua(6)\oplus\bbR$\\
\hline
$\mathfrak{e}_6\oplus\bbR$\\
\hline
\end{tabular}~~,
\end{center}
then the Lie groups $G_k$ corresponding to {\it these} Lie algebras define
the {\it irreducible} Riemannian symmetric spaces $\mathcal{G}_k/G_k$.
Similarly if we replace the third column in the magic square by
\begin{center}
\begin{tabular}{|c|}
\hline
$\spa(3)\oplus\sua(2)$\\
\hline
$\sua(6)\oplus\sua(2)$\\
\hline
$\soa(12)\oplus\sua(2)$\\
\hline
$\mathfrak{e}_7\oplus\sua(2)$\\
\hline
\end{tabular}~~,
\end{center}
then the Lie groups $\mathcal{G}_k$ corresponding to {\it these} Lie
algebras define the {\it irreducible} Riemannian symmetric spaces
$\tilde{\mathcal{G}}_k/\mathcal{G}_k$. Thus starting from the second
and the third column of the table in Theorem \ref{torless}, via the
magic square, we arrived at 12 symmetric spaces
\begin{center}
\begin{tabular}{|c|c|c|}
\hline
$\sug(3)/\sog(3)$&$\spg(3)/\ug(3)$&${\bf F}_4/(\spg(3)\times\sug(2))$\\
\hline
$\sug(3)$&$\sug(6)/{\bf S}(\ug(3)\times\ug(3))$&${\bf E}_6/(\sug(6)\times\sug(2)$\\
\hline
$\sug(6)/\spg(3)$&$\sog(12)/\ug(6)$&${\bf E}_7/(\sog(12)\times\sug(2))$\\
\hline
${\bf E}_6/{\bf F}_4$&${\bf E}_7/({\bf E}_6\times\sog(2))$&${\bf E}_8/({\bf E}_7\times\sug(2))$\\
\hline
\end{tabular}~~.
\end{center}
These 12 symmetric spaces can be considered torsionless models for
special geometries on Riemannian manifolds $M$ with the following
dimensions and structure groups:
\begin{center}
\begin{tabular}{||c|c||c|c||c|c||}
\hline
 $\dim M$&Structure group&$\dim M$&Structure group&$\dim M$&Structure group\\
 $n_k$&$H_k$&$2(n_k+1)$&Extended $G_k$&$4(n_k+2)$&Extended
 $\mathcal{G}_k$\\
\hline\hline
5&$\sog(3)$&12&$\ug(3)$&28&$\spg(3)\times\sug(2)$\\\hline
8&$\sug(3)$&18&${\bf
S}(\ug(3)\times\ug(3))$&40&$\sug(6)\times\sug(2)$\\\hline
14&$\spg(3)$&30&$\ug(6)$&64&$\sog(12)\times\sug(2)$\\\hline 26&${\bf
F}_4$&54&${\bf E}_6\times\sog(2)$&112&${\bf E}_7
\times\sug(2)$\\\hline
\end{tabular}
\end{center}
Quick look at Cartan's list of the irreducible symmetric spaces of
compact type suggests that the special Riemannian geometries
appearing in this list should be supplemented by the two
`exceptional' possibilities: \begin{itemize}
\item[1)] $\dim M=32$, with the structure group $\sog(10)\times\sog(2)$
and with the torsionless model of compact type $M={\bf
E}_6/(\sog(10)\times\sog(2))$
\item[2)] $\dim M=8$, with the structure
group $\sug(2)\times\sug(2)$ and with the torsionless model of
compact type $M={\bf G}_2/(\sug(2)\times\sug(2))$ .
\end{itemize}
Although these two possibilities are not implied by the magic
square, we are convinced that their place is in the above table:
item 1) should stay in the second column for $\dim M$ in the row
between dimensions 30 and 54, and item 2) should stay in the third
column for dim $M$ in the `zeroth' row, before dimension 28.

It is interesting if all these geometries admit characteristic
connection. Also, we do not know what objects in $\bbR^{\dim M}$
reduce the orthogonal groups $\sog(\dim M)$ to the above mentioned
structure groups. Are these symmetric tensors, as it was in the case
of the groups $H_k$?.
\section{Examples in dimension 8}
In the following sections we will briefly discuss the two different 8-dimensional cases, namely: the restricted nearly integrable $\sug(3)$ geometries and the $\sug(2)\times\sug(2)$ geometries. In particular, we provide nontrivial 
examples of 
restricted nearly integrable $\sug(3)$ structures. We also  
explain how to define an $\sug(2)\times\sug(2)$ structure by means of a symmetric tensor of the {\it sixth} order.
\subsection{$\sug(3)$ structures}
It is interesting to note that in the decomposition of
$\bgw^3\bbR^8$ onto the  $\sug(3)$-invariant components (see
Proposition \ref{clastor}) there exists a
1-dimensional subspace ${^8\bgs}^2_1$. This
space, in the adapted coframe of Definition \ref{df:adapted}, is
spanned by a 3-form \be\psi=\tau_1\dz\theta^6+\tau_2\dz\theta^7
+\tau_3\dz\theta^8+\theta^6\dz\theta^7\dz\theta^8,\label{niezmf}\ee
where $(\tau_1,\tau_2,\tau_3)$ are 2-forms
\begin{eqnarray*}
\tau_1&=&\theta^1\dz\theta^4+\theta^2\dz\theta^3+\sqrt{3}\theta^1\dz\theta^5\\
\tau_2&=&\theta^1\dz\theta^3+\theta^4\dz\theta^2+\sqrt{3}\theta^2\dz\theta^5\\
\tau_3&=&\theta^1\dz\theta^2+2\theta^4\dz\theta^3
\end{eqnarray*}
spanning the 3-dimensional irreducible representation
${^5\bgw}^2_3\simeq\soa(3)$ of $\sog(3)$.

Note that the 3-form $\psi$ can be considered in $\bbR^8$ without
any reference to tensor $\ten$. It is remarkable that this 3-form
{\it alone} reduces the ${\bf GL}(8,\bbR)$ to the irreducible
$\sug(3)$ in the same way as $\ten$ does.\footnote{We note that
  $\psi$ is a stable form in dimension 8 and, as such, was considered
  by Nigel Hitchin in \cite{hitchin}} If one thinks that formula
(\ref{niezmf}) is written in the orthonormal coframe $\theta$ then
one gets the reduction from ${\bf GL}(8,\bbR)$ via $\sog(8)$ to the
irreducible $\sug(3)$. Thus, in dimension $8$, the $H_k$ structure
can be defined either in terms of the {\it totally symmetric} $\ten$
or in terms of the {\it totally skew symmetric}
$\psi$.\footnote{Simon Chiossi \cite{chiossi} asks if I there is another 
situation in which a subgroup $H\subset\sog(n)$ of ${\bf GL}(n,\bbR)$ is a
stabilizer of either a totally symmetric tensor or,
independently, of a totally skew symmetric tensor. The answer is
yes. The fifth exterior power of the 14-dimensional
representation ${^{14}\bgs}^2_{14}$ of $\spg(3)$ has a
1-dimensional $\spg(3)$-invariant subspace. Thus, it defines a 5-form 
whose stablizer under the action of ${\bf GL}(14,\bbR)$ includes 
$\spg(3)\subset\sog(14)$. It turns out that this
5-form alone, independently of tensor $\ten$ of
Section \ref{rreppee}, reduces the ${\bf GL}(n,\bbR)$, via $\sog(14)$ to
$\spg(3)$. The explicit expression for this form in the adapted
coframe of Section \ref{hkst} is given in the Appendix B.}
\begin{remark}
In this sense the 3-form $\psi$ and the 2-forms
$(\tau_1,\tau_2,\tau_3)$ play the same role in the relations between
$\sug(3)$ structures in dimension {\it eight} and $\sog(3)$
structures in dimension {\it five} as the 3-form
$$\phi=\sigma_1\dz\theta^5+\sigma_2\dz\theta^6
+\sigma_3\dz\theta^7+\theta^5\dz\theta^6\dz\theta^7$$ and the
self-dual 2-forms
\begin{eqnarray*}
\sigma_1&=&\theta^1\dz\theta^3+\theta^4\dz\theta^2\\
\sigma_2&=&\theta^4\dz\theta^1+\theta^3\dz\theta^2\\
\sigma_3&=&\theta^1\dz\theta^2+\theta^3\dz\theta^4
\end{eqnarray*}
play in the relations \cite{sala} between ${\bf G}_2$ structures in
dimension {\it seven} and $\sug(2)$ structures in dimension {\it
four}.
\end{remark}

We also note that the $\sug(3)$ invariant 3-form $\psi$ can be used
to find the explicit decomposition of an arbitrary 3-form $\omega$
in $\bgw^3\bbR^8$ onto the irreducible components mentioned in
Proposition \ref{clastor}. Indeed, given an $\sug(3)$ structure
$(M,g,\ten)$ on an 8-manifold $M$ we may write an arbitrary 3-form
$\omega$ in the adapted coframe $\theta$ of Definition
\ref{df:adapted} as
$\omega=\tfrac16\omega_{IJK}\theta^I\dz\theta^J\dz\theta^K$. Using
$\psi=\tfrac16\psi_{IJK}\theta^I\dz\theta^J\dz\theta^K$, we
associate with $\omega$ a tensor
$\psi(\omega)_{IJ}=\psi_{IKL}\omega_{JKL}$. Since $\psi(\omega)$ is
an element of $\bigotimes^2\bbR^8$, it may be analyzed by means of
the endomorphism $\hat{\ten}$ naturally associated with
$\ten=\ten_{IJK}\theta^I\theta^J\theta^K$ via (\ref{hatten}). It
follows that the 3-form $\omega$ is
\begin{itemize}
\item in \hspace{.5cm}${^8\bgs^2_{1}}$ \hspace{.5cm}iff
\hspace{.5cm}$\hat{\ten}(\psi(\omega))=20\psi(\omega)$,
\item in \hspace{.5cm}${^8\bgs^2_{8}}$ \hspace{.5cm}iff
\hspace{.5cm}$\hat{\ten}(\psi(\omega))=-6\psi(\omega)$,
\item in \hspace{.5cm}${^8\bgw^2_{20}}$ \hspace{.5cm}iff
\hspace{.5cm}$\hat{\ten}(\psi(\omega))=-8\psi(\omega)$,
\item in \hspace{.5cm}${^8\bgs^2_{27}}$ \hspace{.4cm}iff
\hspace{.5cm}$\hat{\ten}(\psi(\omega))=4\psi(\omega)$.
\end{itemize}
Now, if we have a nearly integrable $\sug(3)$ structure in dimension
$8$, it is easy to check what is the type of its totally skew
symmetric torsion $T_{IJK}$. For this it suffices to consider a
3-form $T=\tfrac16 T_{IJK}\theta^I\dz\theta^J\dz\theta^K$, to
associate with it $\psi(T)$ and to apply the endomorphism
$\hat{\ten}$.
\subsubsection{Structures with $(8+l)$-dimensional symmetry
group}\label{six} In the following we will apply the following
construction.

Let $G$ be an $(8+l)$-dimensional Lie group with a Lie subgroup $H$
of dimension $l\leq 8=\dim\sug(3)$. We arrange a labeling of a left
invariant coframe $(\theta^I,\gamma^\alpha)$ on $G$  in such a way
that the vector fields $X_\alpha$, $\alpha=1,2,\dots,l$, of the
frame $(e_I,X_\alpha)$ dual to $(\theta^I,\gamma^\alpha)$ generate
$H$. Suppose now that the groups $G$ and $H$ are such that the
following tensors
$$g=g_{IJ}\theta^I\theta^J=(\theta^1)^2+(\theta^2)^2+(\theta^3)^2+
(\theta^4)^2+(\theta^5)^2+(\theta^6)^2+(\theta^7)^2+(\theta^8)^2$$
\be\ten=\ten_{IJK}=\tfrac12 \det
\begin{pmatrix}
\theta^5-\sqrt{3}\theta^4&\sqrt{3}(\theta^3+i\theta^8)&\sqrt{3}(\theta^2+i\theta^7)\\
\sqrt{3}(\theta^3-i\theta^8)&\theta^5+\sqrt{3}\theta^4&\sqrt{3}(\theta^1+i\theta^6)\\
\sqrt{3}(\theta^2-i\theta^7)&\sqrt{3}(\theta^1-i\theta^6)&-2\theta^5
\end{pmatrix}\label{naw}\ee
are {\it invariant} on $G$ when Lie transported along the flows of
vector fields $X_\alpha$. In such a case both $g$ and $\ten$ project
to a well defined nondegenerate tensors on the $8$-dimensional
homogeneous space $M=G/H$. These projected tensors, which we will
also respectively denote by $g$ and $\ten$, define the $\sug(3)$
structure $(M,g,\ten)$ on $M$.
\subsubsection{Restricted nearly integrable structures with maximal
symmetry groups} It follows that the {\it restricted} nearly
integrable $\sug(3)$ structures with {\it maximal} symmetry groups
are locally equivalent to the torsionless models of Theorem
\ref{torless}. Thus the possible maximal symmetry groups $G$ are:
$$G_{\lambda>0}=\sug(3)\times\sug(3),\quad
G_{\lambda=0}=\sug(3)\times_\rho\mathbb{R}^8\quad{\rm  or}\quad
G_{\lambda<0}={\bf SL}(3,\bbC).$$ The three cases are
distinguishable by means of the sign of a real constant $\lambda$,
which is related to the Ricci scalar of the
Levi-Civita/characteristic connection of the corresponding
torsionless $\sug(3)$ structure.

We illustrate this statement by using the left invariant coframe
$(\theta^I,\gamma^\alpha)$ on $G$ discussed in the preceding
section. Here it satisfies the following differential system
$$\der\theta^I+\Gamma^I_{~J}\dz\theta^J=0,\quad\quad
\der\Gamma^I_{~J}+\Gamma^I_{~K}\dz\Gamma^K_{~J}=R^I_{~J},
$$
where the characteristic connection matrix $\Gamma^I_J$ is related
to the 1-forms $\gamma^\alpha$, $\alpha=1,2,\dots,8$, via: \be
\Gamma=\gamma^\alpha E_\alpha=\left(\begin{smallmatrix}
0&-\gamma^1&-\gamma^2&-\gamma^3&-\scriptscriptstyle{\sqrt{3}}\gamma^3&-\gamma^4&-\gamma^5&-\gamma^6\\
\gamma^1&0&-\gamma^3&\gamma^2&-\scriptscriptstyle{\sqrt{3}}\gamma^2&-\gamma^5&-\gamma^4-\gamma^7&-\tfrac{\gamma^8}{\scriptscriptstyle{\sqrt{3}}}\\
\gamma^2&\gamma^3&0&2\gamma^1&0&\gamma^6&-\tfrac{\gamma^8}{\scriptscriptstyle{\sqrt{3}}}&-\gamma^7\\
\gamma^3&-\gamma^2&-2\gamma^1&0&0&-\tfrac{\gamma^8}{\scriptscriptstyle{\sqrt{3}}}&-\gamma^6&-2\gamma^5\\
\scriptscriptstyle{\sqrt{3}}\gamma^3&\scriptscriptstyle{\sqrt{3}}\gamma^2&0&0&0&-\gamma^8&\scriptscriptstyle{\sqrt{3}}\gamma^6&0\\
\gamma^4&\gamma^5&-\gamma^6&\tfrac{\gamma^8}{\scriptscriptstyle{\sqrt{3}}}&\gamma^8&0&-\gamma^1&\gamma^2\\
\gamma^5&\gamma^4+\gamma^7&\tfrac{\gamma^8}{\scriptscriptstyle{\sqrt{3}}}&\gamma^6&-\scriptscriptstyle{\sqrt{3}}\gamma^6&\gamma^1&0&-\gamma^3\\
\gamma^6&\tfrac{\gamma^8}{\scriptscriptstyle{\sqrt{3}}}&\gamma^7&2\gamma^5&0&-\gamma^2&\gamma^3&0
\end{smallmatrix}\right),
\label{ko8} \ee and the curvature $R^I_J$ is given by
$$R=-\frac{\lambda}{12}\kappa^\alpha E_\alpha,$$ where the 2-forms
$\kappa^\alpha$ are given by
\begin{eqnarray}
\kappa^1&=&\theta^1\dz\theta^2-2
\theta^3\dz\theta^4+\theta^6\dz\theta^7\nonumber\\
\kappa^2&=&\theta^1\dz\theta^3-
\theta^2\dz\theta^4+\sqrt{3}\theta^2\dz\theta^5-\theta^6\dz\theta^8\nonumber\\
\kappa^3&=&\theta^1\dz\theta^4+\sqrt{3}
\theta^1\dz\theta^5+\theta^2\dz\theta^3+\theta^7\dz\theta^8\nonumber\\
\kappa^4&=&4\theta^1\dz\theta^6+2
\theta^2\dz\theta^7-2\theta^3\dz\theta^8\label{sta}\\
\kappa^5&=&\theta^1\dz\theta^7+
\theta^2\dz\theta^6+2\theta^4\dz\theta^8\nonumber\\
\kappa^6&=&\theta^1\dz\theta^8-
\theta^3\dz\theta^6+\theta^4\dz\theta^7-\sqrt{3}\theta^5\dz\theta^7\nonumber\\
\kappa^7&=&-2\theta^1\dz\theta^6+2
\theta^2\dz\theta^7+4\theta^3\dz\theta^8\nonumber\\
\kappa^8&=&\sqrt{3}(\theta^2\dz\theta^8+
\theta^3\dz\theta^7+\theta^4\dz\theta^6+\sqrt{3}\theta^5\dz\theta^6).\nonumber
\end{eqnarray}
Now, the system guarantees that the Lie derivatives of the
structural tensors $g$ and $\ten$ of (\ref{naw}) with respect to all 
vector fields $X_\alpha$
dual on $G$ to $\gamma^\alpha$ vanish:
$$L_{X_\alpha}g=0,\quad\quad L_{X_\alpha}\ten\equiv
0.$$
Thus, the quotient spaces are equipped with $\sug(3)$
structures locally equivalent to the natural $\sug(3)$ structures on
$$M_{\lambda>0}=\sug(3),\quad M_{\lambda=0}=\bbR^8,\quad
M_{\lambda<0}=\slg(3,\bbC)/\sug(3).$$ Since in such cases the
Riemannian holonomy group is reduced to $\sug(3)$, the Riemannian
curvature coincides with the curvature of the characteristic
connection. The Ricci tensor of these curvatures is Einstein,
$Ric=\lambda g$.
\subsubsection{Examples with 11-dimensional symmetry group and torsion in ${^8\bgs}^2_{27}$} Below we
present a 2-parameter family of restricted nearly integrable
$\sug(3)$ structures which are a local deformation of the
torsionless model $M_{\lambda>0}=\sug(3)$.  These structures have an
11-dimensional symmetry group $G$ described by the Maurer-Cartan
coframe $(\theta^I,\gamma^1,\gamma^2,\gamma^3)$ defined below.
\begin{eqnarray*}
\der\theta^1&=&\gamma^1\dz\theta^2+\gamma^2\dz\theta^3+\gamma^3\dz\theta^4+
\sqrt{3}\gamma^3\dz\theta^5+\frac{k}{\sqrt{3}}\kappa^8\\
\der\theta^2&=&-\gamma^1\dz\theta^1-\gamma^2\dz\theta^4+\sqrt{3}\gamma^2\dz\theta^5+
\gamma^3\dz\theta^3-k\kappa^6\\
\der\theta^3&=&-2\gamma^1\dz\theta^4 -\gamma^2\dz\theta^1-
\gamma^3\dz\theta^2-k\kappa^5 \\
\der\theta^4&=&2\gamma^1\dz\theta^3 +\gamma^2\dz\theta^2-
\gamma^3\dz\theta^1+\frac{k}{2}\kappa^7 \\
\tfrac{1}{\sqrt{3}}\der\theta^5&=&-\gamma^2\dz\theta^2
-\gamma^3\dz\theta^1-\frac{k}{6}(2\kappa^4+\kappa^7)\\
\der\theta^6&=&\gamma^1\dz\theta^7-\gamma^2\dz\theta^8+(2k-t)\kappa^3-2(k+7t)\theta^7\dz\theta^8\\
\der\theta^7&=&-\gamma^1\dz\theta^6+\gamma^3\dz\theta^8+(2k-t)\kappa^2+2(k+7t)\theta^6\dz\theta^8\\
\der\theta^8&=&\gamma^2\dz\theta^6-\gamma^3\dz\theta^7+(2k-t)\kappa^1-2(k+7t)\theta^6\dz\theta^7\\
\der\gamma^1&=&\gamma^2\dz\gamma^3+(k+15t)(t-2k)\kappa^1+(k+15t)(k-t)\theta^6\dz\theta^7\\
\der\gamma^2&=&-\gamma^1\dz\gamma^3+(k+15t)(t-2k)\kappa^2-(k+15t)(k-t)\theta^6\dz\theta^8\\
\der\gamma^3&=&\gamma^1\dz\gamma^2+(k+15t)(t-2k)\kappa^3+(k+15t)(k-t)\theta^7\dz\theta^8
\end{eqnarray*}
The 2-forms $\kappa^\alpha$ appearing here are given in (\ref{sta});
$k$ and $t$ are real constants.

It is easy to check that in all directions spanned by the three vector fields $X_\alpha$ dual
to $\gamma^\alpha$ we have
$$L_{X_\alpha}g=0,\quad\quad L_{X_\alpha}\ten=0,$$
where the structural tensors $g$ and $\ten$ are given by (\ref{naw}). Thus the quotient 8-manifold $M=G/H$, where $H$
is generated by $X_\alpha$, is equipped with an $\sug(3)$ structure
$(M,g,\ten)$. As mentioned above, this structure is restricted
nearly integrable. It has the characteristic connection $\Gamma$
given by (\ref{ko8}) with
$$
\gamma^4=(k-t)(\theta^4+\sqrt{3}\theta^5),\quad\gamma^5=(k-t)\theta^3,
\quad\gamma^6=(k-t)\theta^2,$$
$$\gamma^7=2(t-k)\theta^4,\quad\gamma^8=\sqrt{3}(t-k)\theta^1.
$$
The torsion $T$ of $\Gamma$ is of a pure type. It lies in the
27-dimensional representation ${^8\bgs}^2_{27}$. The torsion 3-form
$T$ reads:
\begin{eqnarray*}
&T=t(\theta^1\dz\theta^2\dz\theta^8+
\theta^1\dz\theta^3\dz\theta^7+\theta^1\dz\theta^4\dz\theta^6+
\sqrt{3}\theta^1\dz\theta^5\dz\theta^6+\theta^2\dz\theta^3\dz\theta^6-\\
&\theta^2\dz\theta^4\dz\theta^7+
\sqrt{3}\theta^2\dz\theta^5\dz\theta^7-2
\theta^3\dz\theta^4\dz\theta^8-15 \theta^6\dz\theta^7\dz\theta^8).
\end{eqnarray*}
Remarkably this form is coclosed, so the Ricci tensor $Ric^\Gamma$
of the characteristic connection is {\it symmetric}. Moreover, it is
diagonal,
$$Ric^\Gamma={\rm
diag}(\lambda,\lambda,\lambda,\lambda,\lambda,\mu,\mu,\mu),
$$
with two constant eigenvalues
$$\lambda=12(k^2+15kt-8t^2),\quad\quad\mu=12(k^2+\tfrac53 kt).$$
These two eigenvalues coincide when $t=0$ and $t=\tfrac53 k$. In the
first case the $\sug(3)$ structure is locally equivalent to the
torsionless model $M_{\lambda>0}=\sug(3)$. The case $$t=\tfrac 53
k$$ is interesting since it provides an example of a restricted
nearly integrable $\sug(3)$ structure with the characteristic
connection $\Gamma$ satisfying the Einstein equations
$$Ric^\Gamma=\tfrac{136}{3}k^2 g$$
and having torsion of a nontrivial pure type ${^8\bgs}^2_{27}$.

We further note that for all values of $t$ and $k$ the Ricci tensor
for the Levi-Civita connection of this structure is also diagonal,
$$Ric^{LC}={\rm
diag}(\lambda',\lambda',\lambda',\lambda',\lambda',\mu',\mu',\mu'),
$$
with eigenvalues given by $\lambda'=\lambda+3t^2$,
$\mu'=\mu+115t^2$. This Ricci tensor is Einstein in the torsionless
case $t=0$ and when $$t=\tfrac{10}{13}k.$$ In this later case the
Ricci tensor reads $$Ric^{LC}=\tfrac{16128}{169}k^2g.$$
\subsubsection{Examples with 9-dimensional symmetry group and vectorial torsion}
Let $G$ be a 9-dimensional group with a left invariant coframe
$(\theta^I,\gamma^1)$ on it as in Section \ref{six}. Let $(e_I,X_1)$
be a basis of vector fields on $G$ dual to $(\theta^I,\gamma^1)$. We
assume that $(\theta^I,\gamma^1)$ satisfies the following
differential system
\begin{eqnarray*}
\der\theta^1&=&\gamma^1\dz\theta^2 - \tfrac12 t_1\theta^1\dz\theta^3
+\tfrac12 t_1\theta^2\dz\theta^4 - \tfrac{1}{2\sqrt{3}} t_1
\theta^2\dz\theta^5 + \tfrac12 t_2\theta^3\dz\theta^7 + \tfrac12 t_2
\theta^4\dz\theta^6\\
\der\theta^2&=&-\gamma^1\dz\theta^1 + \tfrac12
t_1\theta^1\dz\theta^4+\tfrac{1}{2\sqrt{3}} t_1 \theta^1\dz\theta^5+
\tfrac12 t_1\theta^2\dz\theta^3 +\tfrac12 t_2 \theta^3\dz\theta^6 -
\tfrac12 t_2
\theta^4\dz\theta^7\\
\der\theta^3&=&-2\gamma^1\dz\theta^4 + \tfrac{1}{\sqrt{3}}t_1
\theta^4\dz\theta^5\\
\der\theta^4&=&2 \gamma^1\dz\theta^3 - \tfrac{1}{\sqrt{3}}t_1
\theta^3\dz\theta^5\\
\der\theta^5&=&\tfrac{1}{\sqrt{3}} t_1 \theta^1\dz\theta^2 -
\tfrac{1}{\sqrt{3}}t_2 \theta^1\dz\theta^6 -
\tfrac{1}{\sqrt{3}}t_2\theta^2\dz\theta^7 + \tfrac{1}{\sqrt{3}}t_1
\theta^3\dz\theta^4 + \tfrac{1}{\sqrt{3}} t_1 \theta^6\dz\theta^7\\
\der\theta^6&=&\gamma^1\dz\theta^7 + \tfrac12 t_2
\theta^1\dz\theta^4 + \tfrac12 t_2 \theta^2\dz\theta^3 +\tfrac12 t_1
\theta^3\dz\theta^6 -\tfrac12 t_1 \theta^4\dz\theta^7 +
\tfrac{1}{2\sqrt{3}} t_1 \theta^5\dz\theta^7\\
\der\theta^7&=&-\gamma^1\dz\theta^6 + \tfrac12 t_2
\theta^1\dz\theta^3 - \tfrac12 t_2 \theta^2\dz\theta^4 - \tfrac12
t_1\theta^3\dz\theta^7 -\tfrac12 t_1 \theta^4\dz\theta^6
-\tfrac{1}{2 \sqrt{3}} t_1
\theta^5\dz\theta^6\\
\der\theta^8&=& -t_2 \theta^3\dz\theta^4\\
\der\gamma^1&=&-\tfrac16 t_1^2\theta^1\dz\theta^2+\tfrac16
t_1t_2\theta^1\dz\theta^6+\tfrac16
t_1t_2\theta^2\dz\theta^7+\tfrac16(3t_2^2-4t_1^2)\theta^3\dz\theta^4-\tfrac16
t_1^2\theta^6\dz\theta^7.
\end{eqnarray*}
Here the real parameters $t_1, t_2$ are constants.

Let $H$ be a 1-parameter subgroup of $G$ generated by the vector
field $X_1$. Now we consider $g$ and $\ten$ of (\ref{naw}). It is
easy to check that the above differential equations for the system
$(\theta^I,\gamma^1)$ guarantees that on $G$ the Lie derivatives
with respect to $X_1$ of $g$ and $\ten$ identically vanish:
$$L_{X_1}g\equiv 0,\quad\quad\quad L_{X_1}\ten\equiv 0.$$
Thus on the homogeneous space $M=G/H$ we have an $\sug(3)$ structure
$(M,g,\ten)$. This $\sug(3)$ structure has the following properties.
\begin{itemize}
\item It is a restricted nearly integrable structure.
\item It has a 9-dimensional symmetry group $G$.
\item Its characteristic connection is given by (\ref{ko8}),
where
$$
\gamma^2=-\tfrac12 t_1\theta^1,\quad\gamma^3=\tfrac12
t_1\theta^2,\quad\gamma^4=\tfrac12
t_2\theta^4-\tfrac{1}{\sqrt{3}}t_2\theta^5+\tfrac12 t_1\theta^8,
$$
$$
\gamma^5=\tfrac12 t_2\theta^3,\quad\gamma^6=-\tfrac12
t_1\theta^6,\quad\gamma^7=-
t_2\theta^4,\quad\gamma^8=-\tfrac{\sqrt{3}}{2} t_1\theta^7.
$$
\item The skew symmetric torsion $T$ of $\Gamma$ is of
purely `vectorial' type: $T\in{^8\bgs}^2_8$. Explicitly the torsion
3-form is
\begin{eqnarray*}
&T=t_1(-\tfrac{2}{\sqrt{3}}\theta^1\dz\theta^2\dz\theta^5+
\theta^1\dz\theta^6\dz\theta^8+\theta^2\dz\theta^7\dz\theta^8+
\tfrac{1}{\sqrt{3}}\theta^3\dz\theta^4\dz\theta^5-
\tfrac{2}{\sqrt{3}}\theta^5\dz\theta^6\dz\theta^7)+\\
&t_2(\tfrac{1}{\sqrt{3}}\theta^1\dz\theta^5\dz\theta^6+
\tfrac{1}{\sqrt{3}}\theta^2\dz\theta^5\dz\theta^7+
\theta^3\dz\theta^4\dz\theta^8).
\end{eqnarray*}
\item The Ricci tensor $Ric^\Gamma$ of the characteristic
connection is symmetric:
$$Ric^\Gamma=\begin{pmatrix}
-\tfrac13 (4t_1^2+t_2^2) & 0 & 0 & 0 & 0 & 0 &-\tfrac23t_1t_2 & 0\\
0 & -\tfrac13 (4t_1^2+t_2^2)& 0 & 0 & 0 & \tfrac23t_1t_2& 0 & 0 \\
0 & 0 & -\tfrac73t_1^2 & 0 & 0 & 0 & 0 & 0\\
0 & 0 & 0 &-\tfrac73 t_1^2& 0 & 0 & 0 & 0\\
0 & 0 & 0 & 0 & -t_1^2 & 0 & 0 & 0\\
0 & \tfrac23 t_1t_2 & 0 & 0 & 0 & -\tfrac13 (4t_1^2+t_2^2) & 0 &0\\
-\tfrac23t_1t_2 & 0 & 0 & 0 & 0 & 0 & -\tfrac13 (4t_1^2+t_2^2)&0\\
0 & 0 & 0 & 0 & 0 & 0 & 0 & -t_1^2
\end{pmatrix};$$
hence the torsion 3-form $T$ is coclosed.
\item The Ricci tensor $Ric^{LC}$ of the Levi-Civita connection is
$$Ric^{LC}= \left(\begin{smallmatrix}
-\tfrac16 (t_1^2+t_2^2) & 0 & 0 & 0 & 0 & 0 &-\tfrac43t_1t_2 & 0\\
0 & -\tfrac16 (t_1^2+t_2^2)& 0 & 0 & 0 & \tfrac43t_1t_2& 0 & 0 \\
0 & 0 & \tfrac16 (-13t_1^2+3t_2^2) & 0 & 0 & 0 & 0 & 0\\
0 & 0 & 0 &\tfrac16 (-13t_1^2+3t_2^2)& 0 & 0 & 0 & 0\\
0 & 0 & 0 & 0 &\tfrac16 (3t_1^2+2t_2^2) & 0 & 0 & -\tfrac{1}{2\sqrt{3}}t_1t_2\\
0 & \tfrac43 t_1t_2 & 0 & 0 & 0 & -\tfrac16 (t_1^2+t_2^2) & 0 &0\\
-\tfrac43t_1t_2 & 0 & 0 & 0 & 0 & 0 & -\tfrac16 (t_1^2+t_2^2)&0\\
0 & 0 & 0 & 0 & -\tfrac{1}{2\sqrt{3}}t_1t_2 & 0 & 0 & \tfrac12 t_2^2
\end{smallmatrix}\right)
   $$
\end{itemize}
We note that if $t_1=0$ or $t_2=0$ both the Ricci tensors
$Ric^\Gamma$ and $Ric^{LC}$ in the example above are diagonal. Below
we present another 2-parameter family of examples with this
property.

Now the 9-dimensional group $G$ has the basis of the left invariant
forms $(\theta^I,\gamma^1)$  such that:
\begin{eqnarray*}
\der\theta^1&=&\gamma^1\dz\theta^2 - c \theta^2\dz\theta^8
+(t-4c)(\tfrac12\theta^3\dz\theta^7
+\tfrac12\theta^4\dz\theta^6+\frac{3c\sqrt{3}}{6c-t}\theta^5\dz\theta^6)\\
\der\theta^2&=&-\gamma^1\dz\theta^1 + c \theta^1\dz\theta^8
+(t-4c)(\tfrac12\theta^3\dz\theta^6
-\tfrac12\theta^4\dz\theta^7+\frac{3c\sqrt{3}}{6c-t}\theta^5\dz\theta^7)\\
\der\theta^3&=&-2\gamma^1\dz\theta^4 + 2c(\theta^1\dz\theta^7+\theta^2\dz\theta^6+\theta^4\dz\theta^8)\\
\der\theta^4&=&2\gamma^1\dz\theta^3 + 2c(\theta^1\dz\theta^6-\theta^2\dz\theta^7-\theta^3\dz\theta^8)\\
\der\theta^5&=&\frac{6c-t}{\sqrt{3}} (\theta^1\dz\theta^6 +
\theta^2\dz\theta^7)\\
\der\theta^6&=&\gamma^1\dz\theta^7
+(t-4c)(\tfrac12\theta^1\dz\theta^4
+\frac{3c\sqrt{3}}{6c-t}\theta^1\dz\theta^5+\tfrac12\theta^2\dz\theta^3)-
c \theta^7\dz\theta^8\\
\der\theta^7&=&-\gamma^1\dz\theta^6
+(t-4c)(\tfrac12\theta^1\dz\theta^3
-\tfrac12\theta^2\dz\theta^4+\frac{3c\sqrt{3}}{6c-t}\theta^2\dz\theta^5)+
c \theta^6\dz\theta^8\\
\der\theta^8&=& -2c \theta^1\dz\theta^2+(4c-t) \theta^3\dz\theta^4-2c \theta^6\dz\theta^7\\
\der\gamma^1&=&(2c-t)(-c\theta^1\dz\theta^2+\tfrac12 (4c-t)
\theta^3\dz\theta^4-c\theta^6\dz\theta^7),
\end{eqnarray*}
where $c$ and $t$ are constants. The homogeneous space $M=G/H$,
where $H$ is a 1-dimensional subgroup of $G$ generated by $X_1$ dual
to $\gamma^1$, is equipped with a restricted nearly integrable
$\sug(3)$ structure $(M,g,\ten)$ via (\ref{naw}). The characteristic
connection is given by (\ref{ko8}) with
$$
\gamma^2=c\theta^7,\quad\gamma^3=c\theta^6,\quad\gamma^4=\tfrac12
(t-2c)\theta^4+\tfrac{1}{\sqrt{3}}\frac{t^2-18c^2}{6c-t}\theta^5,
$$
$$
\gamma^5=\tfrac12(t-2c)\theta^3,\quad\gamma^6=-c\theta^2,\quad\gamma^7=(2c-t)\theta^4,\quad\gamma^8=\sqrt{3}c\theta^1.
$$
In this 2-parameter family of examples the skew symmetric torsion is
again of purely vectorial type $T\in{^8\bgs}^2_8$; its corresponding
3-form is given by
$$T=t(\tfrac{1}{\sqrt{3}}\theta^1\dz\theta^5\dz\theta^6+
\tfrac{1}{\sqrt{3}}\theta^2\dz\theta^5\dz\theta^7+
\theta^3\dz\theta^4\dz\theta^8).$$ As announced above both the Ricci
tensors are now diagonal for all values of $c$ and $t$. Introducing
the eigenvalues
$$\lambda=\tfrac13
[(6c-t)^2-2t^2],\quad \mu=4c(3c-t),$$
$$\lambda'=\tfrac13
[(6c-t)^2-\tfrac32 t^2],\quad
\mu'=\tfrac13(6c-t)^2,\quad\nu'=\tfrac13 [(6c-t)^2+\tfrac12 t^2]$$
we have
$$
Ric^\Gamma={\rm diag}
(\lambda,\lambda,\mu,\mu,\mu,\lambda,\lambda,\mu),\quad
Ric^{LC}={\rm diag}
(\lambda',\lambda',\nu',\nu',\mu',\lambda',\lambda',\nu').
$$
Of course the group $G$ is a symmetry group of this restricted
nearly integrable $\sug(3)$ structure.

We close this section with the following Theorem, whose proof based
on the Bianchi identities, is purely computational.

\begin{theorem}
Let $(M,g,\ten)$ be an arbitrary restricted nearly integrable
$\sug(3)$ structure in dimension eight. Assume that the torsion $T$
of the characteristic connection of this structure is of purely
vectorial type, $T\in{^8\bgs}^2_8$. Then the 3-form $T$
corresponding to the torsion is coclosed
$$\der(*T)\equiv 0.$$
\end{theorem}
This theorem, in particular, implies that the Ricci tensor of the
characteristic connection for such structures is symmetric.

\subsection{$\sug(2)\times\sug(2)$ structures}
In this section we consider $\sug(2)\times\sug(2)$ structures in
dimension eight modeled on the torsionless structure ${\bf
G}_2/(\sug(2)\times\sug(2)$. The approach presented here should be
useful in studies of the other exceptional case concerning with the 
$\sog(10)\times\sog(2)$ structures in dimension 32.

In full analogy with $H_k$ structures we start with the
identification of $\bbR^8$ with a space of the
antisymmetric block matrices $M_{7\times
7}(\mathbb{R})\ni\iota(\vec{X})=\begin{pmatrix} {\bf 0}_{3\times
3}&\alpha\\
-\alpha^t&{\bf 0}_{4\times 4},
\end{pmatrix},$ in which the matrices $\alpha\in M_{3\times 4}(\bbR)$
have 3 rows and 4 columns. The  entries of $\alpha$ satisfy the
following four relations
\begin{eqnarray} &&\alpha_{16}-\alpha_{34}-\alpha_{25}=0,\quad
\alpha_{26}-\alpha_{37}+\alpha_{15}=0,\nonumber\\
&&\alpha_{36}+\alpha_{27}+\alpha_{14}=0,\quad
\alpha_{35}+\alpha_{17}-\alpha_{24}=0\label{relg2}\end{eqnarray}
These four relations reduce the 12 free parameters standing in
an arbitrary $3\times 4$ matrix to 8 parameters. Now defining
\be \bbM^8=\{~\iota(\vec{X})\in M_{7\times
7}(\mathbb{R}):~~\iota(\vec{X})=\begin{pmatrix} {\bf 0}_{3\times
3}&\alpha\\
-\alpha^t&{\bf 0}_{4\times 4},
\end{pmatrix}~{\rm with}~ \alpha~{\rm satisfying}~{\rm (\ref{relg2})}
~\},\label{r5} \ee
we have an isomorphism
$\iota:\bbR^8\to\bbM^8$ between the vector spaces $\bbR^8$ and
$\bbM^8$.

Now we define a representation $\rho$ of the group
$\sug(2)\times\sug(2)$ in $\bbR^8$, which will enable us to define
an $\sug(2)\times\sug(2)$ structure in dimension eight.

We use two different representations of $\sug(2)$ in dimension
seven. The representation $\rho_1$ generated by $7\times 7$
matrices
$$h_{~i}=\exp(t_i s_i),\quad i=1,2,3,\quad t_i\in\bbR\quad {\rm(no~~ summation!),}$$
such that
$$s_1=\left(\begin{smallmatrix}
0&0&0&0&0&0&0\\0&0&0&0&0&0&0\\0&0&0&0&0&0&0\\
0&0&0&0&\tfrac12&0&0\\0&0&0&-\tfrac12&0&0&0\\0&0&0&0&0&0&-\tfrac12\\0&0&0&0&0&\tfrac12&0\end{smallmatrix}\right),
\quad s_2=\left(\begin{smallmatrix}
0&0&0&0&0&0&0\\0&0&0&0&0&0&0\\0&0&0&0&0&0&0\\
0&0&0&0&0&\tfrac12&0\\0&0&0&0&0&0&\tfrac12\\0&0&0&-\tfrac12&0&0&0\\0&0&0&0&-\tfrac12&0&0\end{smallmatrix}\right),
\quad s_3=\left(\begin{smallmatrix}
0&0&0&0&0&0&0\\0&0&0&0&0&0&0\\0&0&0&0&0&0&0\\
0&0&0&0&0&0&\tfrac12\\0&0&0&0&0&-\tfrac12&0\\0&0&0&0&\tfrac12&0&0\\0&0&0&-\tfrac12&0&0&0\end{smallmatrix}\right),$$
and the representation $\rho_2$ generated by $7\times 7$ matrices
$$\chi_i=\exp(\tau_i \sigma_i),\quad i=1,2,3,\quad\tau_i\in\bbR\quad {\rm(no~~ summation!),}$$
such that
$$\sigma_1=\left(\begin{smallmatrix}
0&-1&0&0&0&0&0\\1&0&0&0&0&0&0\\0&0&0&0&0&0&0\\
0&0&0&0&0&0&-\tfrac12\\0&0&0&0&0&-\tfrac12&0\\0&0&0&0&\tfrac12&0&0\\0&0&0&\tfrac12&0&0&0\end{smallmatrix}\right),
\quad \sigma_2=\left(\begin{smallmatrix}
0&0&1&0&0&0&0\\0&0&0&0&0&0&0\\-1&0&0&0&0&0&0\\
0&0&0&0&0&\tfrac12&0\\0&0&0&0&0&0&-\tfrac12\\0&0&0&-\tfrac12&0&0&0\\0&0&0&0&\tfrac12&0&0\end{smallmatrix}\right),
\quad\sigma_3=\left(\begin{smallmatrix}
0&0&0&0&0&0&0\\0&0&1&0&0&0&0\\0&-1&0&0&0&0&0\\
0&0&0&0&-\tfrac12&0&0\\0&0&0&\tfrac12&0&0&0\\0&0&0&0&0&0&-\tfrac12\\0&0&0&0&0&\tfrac12&0\end{smallmatrix}\right).$$
We note that
$$[s_j,s_k]=\epsilon_{ijk}s_i,\quad
[\sigma_j,\sigma_k]=\epsilon_{ijk}\sigma_i,\quad
[s_i,\sigma_j]=0,\quad i,j,k=1,2,3,$$ where $\epsilon_{ijk}$ is the
Levi-Civita symbol in 3-dimensions.

Now, we consider all the $7\times 7$-matrices of the form
$$h=h(t_1,~t_2,~t_3,~\tau_1,~\tau_2~\tau_3)=h_1~h_2~h_3~\chi_1~\chi_2~\chi_3.$$
They constitute a 7-dimensional representation $\rho_7$ of the full
group $\sug(2)\times\sug(2)$.

Remarkably, $h\iota(\vec{X})h^t$ is an element of $\bbM^8$ for all the elements
$\iota(\vec{X})$ of $\bbM^8$. Moreover, due to the fact that $[s_i,\sigma_j]=0$
for all $i,j$, the map
$$(\sug(2)\times\sug(2))\times\bbM^8\ni\quad (h,\iota(\vec{X}))\mapsto h~\iota(\vec{X})~h^t\quad\in\bbM^8$$ is a good action of
$\sug(2)\times\sug(2)$ on $\bbM^8$. Thus, using the isomorphism
$\iota$ we get the 8-dimensional representation $\rho$ of
$\sug(2)\times\sug(2)$ given by
\be\bbR^8\ni\vec{X}\quad\mapsto\quad\rho(h)\vec{X}=\iota^{-1}[~h~\iota
(\vec{X})~h^t~]\quad\in\bbR^8.\label{so35}\ee

Given an element $\vec{X}\in\mathbb{R}^8$ we consider its
characteristic polynomial
\begin{eqnarray}
P_{\vec{X}}(\lambda)&=&{\rm det}(\iota(\vec{X})-\lambda
I)\nonumber\\
&=&-\lambda^7-6g(\vec{X},\vec{X})\lambda^5-9g(\vec{X},\vec{X})^2\lambda^3+2
\gamma(\vec{X},\vec{X},\vec{X},\vec{X},\vec{X},\vec{X})\lambda\label{med}.
\end{eqnarray}
This polynomial is invariant under the $\sug(2)\times\sug(2)$-action
given by the representation $\rho$ of (\ref{so35}),
$$
P_{\rho(h)\vec{X}}(\lambda)=P_{\vec{X}}(\lambda).
$$
Thus, all the coefficients of $P_{\vec{X}}(\lambda)$, which are
multilinear in $\vec{X}$, are $\sug(2)\times\sug(2)$-invariant.

It is convenient to use a basis ${\bf e}_I$ in $\bbR^8$ such that
the isomorphism $\iota:\bbR^8\to\bbM^8$ takes the form:
$$\vec{X}=x^I{\bf e}_I\mapsto \iota(\vec{X})=\left(\begin{smallmatrix}
0& 0&0&-x^3+ \sqrt{3}x^4&-x^5 + \sqrt{3}x^6&-x^7 + \sqrt{3}x^8& -2 x^1\\
0& 0& 0&-x^1 - \sqrt{3}x^2&x^7 + \sqrt{3}x^8& -x^5 - \sqrt{3}x^6& 2 x^3\\
0& 0&0&-2x^7&x^1 - \sqrt{3}x^2&-x^3 - \sqrt{3} x^4& -2 x^5\\
x^3 - \sqrt{3} x^4& x^1 + \sqrt{3} x^2& 2 x^7& 0& 0& 0&0\\
x^5 - \sqrt{3}x^6& -x^7 - \sqrt{3} x^8& -x^1 + \sqrt{3}x^2& 0& 0& 0& 0\\
x^7 - \sqrt{3} x^8& x^5 + \sqrt{3}x^6& x^3 + \sqrt{3} x^4& 0& 0& 0& 0\\
2 x^1& -2 x^3& 2 x^5& 0& 0& 0& 0
\end{smallmatrix}\right).
$$
With this choice the bilinear form $g$ of (\ref{med}) reads:
$$g(\vec{X},\vec{X})=(x^1)^2+(x^2)^2+(x^3)^2+(x^4)^2+(x^5)^2+(x^6)^2+(x^7)^2+(x^8)^2.$$
The 6-linear form $\gamma$ of (\ref{med}) defines a tensor
$\gamma_{IJKLMN}$ via $$
\gamma(\vec{X},\vec{X},\vec{X},\vec{X},\vec{X},\vec{X})=\gamma_{IJKLMN}x^I
x^Jx^Kx^Lx^Mx^N.$$ This tensor has the following
properties.
\begin{itemize}
\item It reduces ${\bf GL}(8,\bbR)$, via $\sog(8)$, to
$\sug(2)\times\sug(2)$.
\item The 6th order polynomial
$$\Phi=\gamma(\vec{X},\vec{X},\vec{X},\vec{X},\vec{X},\vec{X})$$
of variables $x^I$, $I=1,2,\dots,8$, satisfies
\begin{eqnarray*}
&&{\rm a)}\quad\quad\quad  \triangle \Phi=-72~ g(\vec{X},\vec{X})^2\\
&&{\rm b)}\quad\quad\quad |\vec{\nabla} \Phi|^2=-72~\Phi~ g(\vec{X},\vec{X})^2\\
&&{\rm c)}\quad\quad\quad \vec{X}\vec{\nabla}\Phi=6\Phi.
\end{eqnarray*}
\end{itemize}
The properties a)-c) show that $\Phi$ can not be interpreted as the Cartan polynomial (\ref{carmun1})-(\ref{carmun2}) defining an isoparametric hypersurface in $S^7$. But we can modify it so that the redefined polynomial satisfies (\ref{carmun1})-(\ref{carmun2}). Indeed, using properties a)-c) it is easy to see that 
the 6th order homogeneous polynomial
$
F=\Phi+g(\vec{X},\vec{X})^3
$
is a solution of 
\begin{eqnarray*}
&&{\rm Cii)}\quad\quad\quad \triangle F=0\\
&&{\rm Ciii)}\quad\quad\quad |\vec{\nabla}F|^2=6^2~g(\vec{X},\vec{X})^5.
\end{eqnarray*}
Thus, via (\ref{carmun3}), the polynomial $F$ defines an isoparametric hypersurface in $S^7$ which has $p=6$ distinct constant principal eigenvalues. Note that since both $\Phi$ and $g(\vec{X},\vec{X})$ are $\sug(2)\times\sug(2)$ invariant, the polynomial $F$ also is. Hence a stabilizer, under the action of ${\bf GL}(8,\bbR)$, of a 6th order symmetric tensor $\ten_{IJKLMN}$ defined by
\be
F=\Phi+g(\vec{X},\vec{X})^3=\ten(\vec{X},\vec{X},\vec{X},\vec{X},\vec{X},\vec{X})=\ten_{IJKLMN}x^I
x^Jx^Kx^Lx^Mx^N \label{ten6}\ee
contains the group $\sug(2)\times\sug(2)$. Actually we have the following proposition.
\begin{proposition}
The 6th order symmetric tensor $\ten_{IJKLMN}$ defined above reduces the ${\bf GL}(8,\bbR)$ group, via $\sog(8)$, to the irreducible $\sug(2)\times\sug(2)$ associated with the representation $\rho$ of (\ref{so35}).
\end{proposition}
%\subsection{$\spg(3)$ geometries in dimension 14}
%\subsection{${\bf F}_4$ geometries in dimension 26}

Following the case of $H_k$ structures we use the tensor
$\ten_{IJKLMN}$ of (\ref{ten6}) to define an endomorphism
\be\hat{\ten}:{\textstyle \bigotimes}^2\bbR^8
\longrightarrow{\textstyle\bigotimes}^2\bbR^8,\label{hatten6}\ee
$$W^{IK}\stackrel{\hat{\ten}}{\longmapsto} \tfrac{5^2}{2^5}~\ten_{IJMNPQ}\ten_{KLMNPQ}
    W^{JL},$$
which preserves the decomposition
    $\bigotimes^2\bbR^8=\bgw^2\bbR^8\oplus\bgs^2\bbR^8$.
Its eigenspaces, are $\sug(2)\times\sug(2)$ invariant and define
representations of dimension 1,5,6,7,9,15,21. Explicitly we have the
following proposition.
\begin{proposition} \label{pr:su3rep}~\\
The $\sug(2)\times\sug(2)$ irreducible decomposition of
$\bgt^2\bbR^8$ is given by
$$\bgt^2\bbR^8 = \bgs^2_1 \oplus \bgs^2_5 \oplus \bgs^2_9 \oplus\bgs^2_{21} \oplus
\bgw^2_6\oplus\bgw^2_7\oplus\bgw^2_{15},$$ where
\begin{eqnarray*}
  &&\bgs^2_1       = \{~S\in\bgt^2\bbR^8~|~ \hat{\ten} (S)= 175 \cdot S~\} = \{S=\lambda
  \cdot g,\;\lambda\in\bbR~\}, \\
  &&\bgs^2_5       = \{~S\in\bgt^2\bbR^8~|~ \hat{\ten} (S)= -21 \cdot S~\}, \\
 && \bgw^2_6 = \{ ~F\in\bgt^2 \bbR^8~|~  \hat{\ten} (F) = 35 \cdot F~\}~
=~\sua(2)\oplus\sua(2),\\
&& \bgw^2_7 = \{ ~F\in\bgt^2 \bbR^8~|~  \hat{\ten} (F) = -25 \cdot
F~\},\\
 &&\bgs^2_9       = \{~S\in\bgt^2\bbR^8~|~ \hat{\ten} (S)= 7 \cdot S~\}, \\
 && \bgw^2_{15} = \{ ~F\in\bgt^2 \bbR^8~|~  \hat{\ten} (F) = -49 \cdot
F~\},\\
&&\bgs^2_{21}       = \{~S\in\bgt^2\bbR^8~|~ \hat{\ten} (S)= 27
\cdot S~\}.
\end{eqnarray*}
The \emph{real} vector spaces $\bgw^2_i\subset\bgw^2\bbR^8$ and
$\bgs^2_j\subset\bgs^2\bbR^8$ of respective dimensions $i$ and $j$
are irreducible representations of the group $\sug(2)\times\sug(2)$.
\end{proposition}
\begin{remark}
Note that $\bgw^2_6$, is isomorphic to the Lie algebra
$\sua(2)\oplus\sua(2)$ represented as the Lie subalgebra of $8\times
8$ matrices. In this 8-dimensional representation $\rho'$ the basis
of the two $\sua(2)$ algebras are, respectively, $\Sigma^L_i$ and
$\Sigma^R_i$, $i=1,2,3$, where
$$\Sigma^L_1=\tfrac14\left(\begin{smallmatrix}
0&0&-5&\sqrt{3}&0&0&0&0\\0&0&-\sqrt{3}&3&0&0&0&0\\5&\sqrt{3}&0&0&0&0&0&0\\
-\sqrt{3}&-3&0&0&0&0&0&0\\0&0&0&0&0&0&-2&0\\0&0&0&0&0&0&0&6\\0&0&0&0&2&0&0&0\\
0&0&0&0&0&-6&0&0\end{smallmatrix}\right),
\quad\Sigma^L_2=\tfrac14\left(\begin{smallmatrix}
0&0&0&0&-5&\sqrt{3}&0&0\\0&0&0&0&\sqrt{3}&-3&0&0\\0&0&0&0&0&0&-1&-\sqrt{3}\\
0&0&0&0&0&0&3\sqrt{3}&-3\\5&-\sqrt{3}&0&0&0&0&0&0\\-\sqrt{3}&3&0&0&0&0&0&0\\
0&0&1&-3\sqrt{3}&0&0&0&0\\0&0&\sqrt{3}&3&0&0&0&0\end{smallmatrix}\right)$$
$$\Sigma^L_3=\tfrac14\left(\begin{smallmatrix}
0&0&0&0&0&0&-1&\sqrt{3}\\0&0&0&0&0&0&-3\sqrt{3}&-3\\0&0&0&0&5&\sqrt{3}&0&0\\
0&0&0&0&\sqrt{3}&3&0&0\\0&0&-5&-\sqrt{3}&0&0&0&0\\0&0&-\sqrt{3}&-3&0&0&0&0\\1&3\sqrt{3}&0&0&0&0&0&0\\
-\sqrt{3}&3&0&0&0&0&0&0\end{smallmatrix}\right),$$
$$\Sigma^R_1=\tfrac14\left(\begin{smallmatrix}
0&0&0&0&0&0&-1&\sqrt{3}\\0&0&0&0&0&0&\sqrt{3}&1\\0&0&0&0&1&\sqrt{3}&0&0\\
0&0&0&0&\sqrt{3}&-1&0&0\\0&0&-1&-\sqrt{3}&0&0&0&0\\0&0&-\sqrt{3}&1&0&0&0&0\\1&-\sqrt{3}&0&0&0&0&0&0\\
-\sqrt{3}&-1&0&0&0&0&0&0\end{smallmatrix}\right),\quad\Sigma^R_2=\tfrac14\left(\begin{smallmatrix}
0&0&0&0&1&-\sqrt{3}&0&0\\0&0&0&0&-\sqrt{3}&-1&0&0\\0&0&0&0&0&0&1&\sqrt{3}\\
0&0&0&0&0&0&\sqrt{3}&-1\\-1&\sqrt{3}&0&0&0&0&0&0\\\sqrt{3}&1&0&0&0&0&0&0\\
0&0&-1&-\sqrt{3}&0&0&0&0\\0&0&-\sqrt{3}&1&0&0&0&0\end{smallmatrix}\right),$$
$$\Sigma^R_3=\tfrac14\left(\begin{smallmatrix}
0&0&1&-\sqrt{3}&0&0&0&0\\0&0&\sqrt{3}&1&0&0&0&0\\-1&-\sqrt{3}&0&0&0&0&0&0\\
\sqrt{3}&-1&0&0&0&0&0&0\\0&0&0&0&0&0&2&0\\0&0&0&0&0&0&0&2\\0&0&0&0&-2&0&0&0\\
0&0&0&0&0&-2&0&0\end{smallmatrix}\right)$$ and we have
$$[\Sigma^L_j,\Sigma^L_k]=-\epsilon_{ijk}\Sigma^L_i,\quad
[\Sigma^R_j,\Sigma^R_k]=-\epsilon_{ijk}\Sigma^R_i,\quad
[\Sigma^L_i,\Sigma^R_j]=0,\quad i,j,k=1,2,3.$$ Of course the Lie
algebra representation $\rho'$ is the derivative of the
8-dimensional representation $\rho$ of $\sug(2)\times\sug(2)$
considered in (\ref{so35}).
\end{remark}
Now, we define the $\sug(2)\times\sug(2)$ structure on an 8-manifold
as a structure equipped with the Riemannian metric $g$ an the
6-tensor $\ten$, which in an orthonormal coframe $\theta^I$ is given
by
$\ten=\ten_{IJKLMN}\theta^I\theta^J\theta^K\theta^L\theta^M\theta^N$
with $\ten_{IJKLMN}$ of (\ref{ten6}). Having this, we may use the
above proposition (and the basis of $\bgw^2_6=\sua(2)\oplus\sua(2)$
given above) as the starting point for addressing the question about
the description of such structures in terms of
$(\sua(2)\oplus\sua(2))$-valued connections. Regardless of the open
question if and when the characteristic connection for such
structures exists, the torsionless models here will locally be
isometric to the symmetric spaces ${\bf
G}_2/(\sug(2)\times\sug(2))$,
$\bbR^8=[(\sug(2)\times\sug(2))\times_\rho\bbR^8]/(\sug(2)\times\sug(2))$
and ${\bf G}^2_2/(\sug(2)\times\sug(2))$ with the standard
$\sug(2)\times\sug(2)$ structure on them.
\section{Acknowledgements}
This paper owes much to the discussion with Robert Bryant I had
during the RIMS Symposium "Developments of Cartan Geometry and
Related Mathematical Problems" held in Kyoto in October 2005. I am
very grateful to Tohru Morimoto for inviting me to this symposium.

I warmly thank Thomas Friedrich for directing my attention on
special Riemannian geometries, and to Ilka Agricola for answering
many questions concerned with the representation theory.

\renewcommand{\theequation}{A-\arabic{equation}}
\setcounter{equation}{0}
\section*{Appendix A. Basae for the
$n_k$-dimensional representations of the Lie algebras ${\mathfrak
h}_k$} Here we give the explicit formulae for the generic elements
$X_{n_k}$ of the Lie algebras ${\mathfrak h}_k$ in terms of the
$n_k\times n_k$ antisymmetric matrices. We denote the basis of the
Lie algebra ${\mathfrak h}_k$ by $E_\alpha$, $\alpha=1,2,\dots,\dim
H_k$ and write
$$X_{n_k}=x^\alpha E_\alpha.$$
The explicit form of the matrices $E_\alpha$ for each value of
$n_k=5,8,14$ and $26$ can be read off from the formulae below.

For $n_k=5$ we have:
\begin{equation}
X_5=\left(\begin{smallmatrix}
0&-x^1&-x^2&-x^3&-\scriptscriptstyle{\sqrt{3}}x^3\\
x^1&0&-x^3&x^2&-\scriptscriptstyle{\sqrt{3}}x^2\\
x^2&x^3&0&2x^1&0\\
x^3&-x^2&-2x^1&0&0\\
\scriptscriptstyle{\sqrt{3}}x^3&\scriptscriptstyle{\sqrt{3}}x^2&0&0&0
\end{smallmatrix}\right).\label{x5}
\end{equation}

For $n_k=8$ we have:
\begin{equation}\label{x8}
X_8=\left(\begin{smallmatrix}
0&-x^1&-x^2&-x^3&-\scriptscriptstyle{\sqrt{3}}x^3&-x^4&-x^5&-x^6\\
x^1&0&-x^3&x^2&-\scriptscriptstyle{\sqrt{3}}x^2&-x^5&-x^4-x^7&-\tfrac{x^8}{\scriptscriptstyle{\sqrt{3}}}\\
x^2&x^3&0&2x^1&0&x^6&-\tfrac{x^8}{\scriptscriptstyle{\sqrt{3}}}&-x^7\\
x^3&-x^2&-2x^1&0&0&-\tfrac{x^8}{\scriptscriptstyle{\sqrt{3}}}&-x^6&-2x^5\\
\scriptscriptstyle{\sqrt{3}}x^3&\scriptscriptstyle{\sqrt{3}}x^2&0&0&0&-x^8&\scriptscriptstyle{\sqrt{3}}x^6&0\\
x^4&x^5&-x^6&\tfrac{x^8}{\scriptscriptstyle{\sqrt{3}}}&x^8&0&-x^1&x^2\\
x^5&x^4+x^7&\tfrac{x^8}{\scriptscriptstyle{\sqrt{3}}}&x^6&-\scriptscriptstyle{\sqrt{3}}x^6&x^1&0&-x^3\\
x^6&\tfrac{x^8}{\scriptscriptstyle{\sqrt{3}}}&x^7&2x^5&0&-x^2&x^3&0
\end{smallmatrix}\right).
\end{equation}

For $n_k=14$ we have:
\begin{equation}
X_{14}=\label{x14}
\end{equation}
$$
\left(\begin{smallmatrix}
0&-x^1&-x^2&-x^3&-\scriptscriptstyle{\sqrt{3}}x^3&-x^4&-x^5&-x^6&-x^9&-x^{10}&-x^{11}&-x^{12}&-x^{13}&-x^{14}\\
x^1&0&-x^3&x^2&-\scriptscriptstyle{\sqrt{3}}x^2&-x^5&-x^4-x^7&-\tfrac{x^8}{\scriptscriptstyle{\sqrt{3}}}&-x^{11}&-x^{12}&-x^9-x^{15}&-x^{10}-x^{16}&-\tfrac{x^{17}}{\scriptscriptstyle{\sqrt{3}}}&-\tfrac{x^{18}}{\scriptscriptstyle{\sqrt{3}}}\\
x^2&x^3&0&2x^1&0&x^6&-\tfrac{x^8}{\scriptscriptstyle{\sqrt{3}}}&-x^7&x^{13}&x^{14}&-\tfrac{x^{17}}{\scriptscriptstyle{\sqrt{3}}}&-\tfrac{x^{18}}{\scriptscriptstyle{\sqrt{3}}}&-x^{15}&-x^{16}\\
x^3&-x^2&-2x^1&0&0&-\tfrac{x^8}{\scriptscriptstyle{\sqrt{3}}}&-x^6&-2x^5&-\tfrac{x^{17}}{\scriptscriptstyle{\sqrt{3}}}&-\tfrac{x^{18}}{\scriptscriptstyle{\sqrt{3}}}&-x^{13}&-x^{14}&-2x^{11}&-2x^{12}\\
\scriptscriptstyle{\sqrt{3}}x^3&\scriptscriptstyle{\sqrt{3}}x^2&0&0&0&-x^8&\scriptscriptstyle{\sqrt{3}}x^6&0&-x^{17}&-x^{18}&\scriptscriptstyle{\sqrt{3}}x^{13}&\scriptscriptstyle{\sqrt{3}}x^{14}&0&0\\
x^4&x^5&-x^6&\tfrac{x^8}{\scriptscriptstyle{\sqrt{3}}}&x^8&0&-x^1&x^2&-x^{19}&-x^{20}&x^{12}&-x^{11}&-x^{14}&x^{13}\\
x^5&x^4+x^7&\tfrac{x^8}{\scriptscriptstyle{\sqrt{3}}}&x^6&-\scriptscriptstyle{\sqrt{3}}x^6&x^1&0&-x^3&x^{12}&-x^{11}&x^{16}-x^{19}&-x^{15}-x^{20}&-\tfrac{x^{18}}{\scriptscriptstyle{\sqrt{3}}}&\tfrac{x^{17}}{\scriptscriptstyle{\sqrt{3}}}\\
x^6&\tfrac{x^8}{\scriptscriptstyle{\sqrt{3}}}&x^7&2x^5&0&-x^2&x^3&0&-x^{14}&x^{13}&-\tfrac{x^{18}}{\scriptscriptstyle{\sqrt{3}}}&\tfrac{x^{17}}{\scriptscriptstyle{\sqrt{3}}}&z^{10}&-z^9\\
x^9&x^{11}&-x^{13}&\tfrac{x^{17}}{\scriptscriptstyle{\sqrt{3}}}&x^{17}&x^{19}&-x^{12}&x^{14}&0&x^{21}&-x^1&x^5&x^2&-x^6\\
x^{10}&x^{12}&-x^{14}&\tfrac{x^{18}}{\scriptscriptstyle{\sqrt{3}}}&x^{18}&x^{20}&x^{11}&-x^{13}&-x^{21}&0&-x^5&-x^1&x^6&x^2\\
x^{11}&x^9+x^{15}&\tfrac{x^{17}}{\scriptscriptstyle{\sqrt{3}}}&x^{13}&-\scriptscriptstyle{\sqrt{3}}x^{13}&-x^{12}&-x^{16}+x^{19}&\tfrac{x^{18}}{\scriptscriptstyle{\sqrt{3}}}&x^1&x^5&0&x^7+x^{21}&-x^3&-\tfrac{x^{8}}{\scriptscriptstyle{\sqrt{3}}}\\
x^{12}&x^{10}+x^{16}&\tfrac{x^{18}}{\scriptscriptstyle{\sqrt{3}}}&x^{14}&-\scriptscriptstyle{\sqrt{3}}x^{14}&x^{11}&x^{15}+x^{20}&-\tfrac{x^{17}}{\scriptscriptstyle{\sqrt{3}}}&-x^5&x^1&-x^7-x^{21}&0&\tfrac{x^{8}}{\scriptscriptstyle{\sqrt{3}}}&-x^3\\
x^{13}&\tfrac{x^{17}}{\scriptscriptstyle{\sqrt{3}}}&x^{15}&2x^{11}&0&x^{14}&\tfrac{x^{18}}{\scriptscriptstyle{\sqrt{3}}}&-z^{10}&-x^2&-x^6&x^3&-\tfrac{x^{8}}{\scriptscriptstyle{\sqrt{3}}}&0&z^4\\
x^{14}&\tfrac{x^{18}}{\scriptscriptstyle{\sqrt{3}}}&x^{16}&2x^{12}&0&-x^{13}&-\tfrac{x^{17}}{\scriptscriptstyle{\sqrt{3}}}&z^9&x^6&-x^2&\tfrac{x^{8}}{\scriptscriptstyle{\sqrt{3}}}&x^3&-z^4&0
\end{smallmatrix}\right),
$$
where $z^4=x^4+x^7+x^{21}$, $z^9=x^{9}+x^{15}+x^{20}$,
$z^{10}=x^{10}+x^{16}-x^{19}$.

The size of the formula for $X_{14}$ forces us to skip the
$26$-dimensional representation of $\mathfrak{f}_4$. It can be
easily obtained by looking for the Lie algebra element stabilizing
$\ten^1$ of Proposition \ref{pb}.

\renewcommand{\theequation}{B-\arabic{equation}}
\setcounter{equation}{0}
\section*{Appendix B. A 5-form reducing ${\bf GL}(14,\bbR)$ to
  $\spg(3)\subset\sog(14)$}
An explicit expression for a non-zero element $\phi$ of the
only 1-dimensional $\spg(3)$-invariant subspace in
$\bgw^5~{^{14}\bgs}^2_{14}$ is written below. This 5-form reduces 
${\bf GL}(14,\bbR)$ to  $\spg(3)\subset\sog(3)$ and, 
in the adpated coframe of Section \ref{hkst}, contains 129 terms as follows:
\begin{eqnarray*}&&\phi=\\
&& 120\ \sqrt{3}\ \theta^1\dz\theta^2\dz\theta^3\dz\theta^4\dz\theta^5 - 
    240\ \sqrt{3}\ \theta^1\dz\theta^2\dz\theta^5\dz\theta^6\dz\theta^7 - 
    192\ \sqrt{3}\
    \theta^1\dz\theta^2\dz\theta^5\dz\theta^9\dz\theta^{11} -\\
&& 
    144\ \sqrt{3}\ \theta^1\dz\theta^2\dz\theta^5\dz\theta^{10}\dz\theta^{12} - 
    72\ \theta^1\dz\theta^2\dz\theta^6\dz\theta^9\dz\theta^{14} + 
    144\ \theta^1\dz\theta^2\dz\theta^6\dz\theta^{10}\dz\theta^{13} - \\&&
    54\ \theta^1\dz\theta^2\dz\theta^7\dz\theta^{11}\dz\theta^{14} + 
    72\ \theta^1\dz\theta^2\dz\theta^7\dz\theta^{12}\dz\theta^{13} - 
    360\ \theta^1\dz\theta^2\dz\theta^8\dz\theta^9\dz\theta^{10} - \\&&
    162\ \theta^1\dz\theta^2\dz\theta^8\dz\theta^{11}\dz\theta^{12} - 
    18\ \theta^1\dz\theta^2\dz\theta^8\dz\theta^{13}\dz\theta^{14} - 
    360\ \theta^1\dz\theta^3\dz\theta^4\dz\theta^6\dz\theta^8 - \\&&
    144\ \theta^1\dz\theta^3\dz\theta^4\dz\theta^9\dz\theta^{13} - 
    72\ \theta^1\dz\theta^3\dz\theta^4\dz\theta^{10}\dz\theta^{14} - 
    120\ \sqrt{3}\ \theta^1\dz\theta^3\dz\theta^5\dz\theta^6\dz\theta^8 -\\&& 
    48\ \sqrt{3}\ \theta^1\dz\theta^3\dz\theta^5\dz\theta^9\dz\theta^{13} - 
    24\ \sqrt{3}\ \theta^1\dz\theta^3\dz\theta^5\dz\theta^{10}\dz\theta^{14} - 
    216\ \theta^1\dz\theta^3\dz\theta^6\dz\theta^9\dz\theta^{12} + \\&&
    288\ \theta^1\dz\theta^3\dz\theta^6\dz\theta^{10}\dz\theta^{11} - 
    360\ \theta^1\dz\theta^3\dz\theta^7\dz\theta^9\dz\theta^{10} - 
    162\ \theta^1\dz\theta^3\dz\theta^7\dz\theta^{11}\dz\theta^{12} - \\&&
    18\ \theta^1\dz\theta^3\dz\theta^7\dz\theta^{13}\dz\theta^{14} - 
    54\ \theta^1\dz\theta^3\dz\theta^8\dz\theta^{11}\dz\theta^{14} + 
    72\ \theta^1\dz\theta^3\dz\theta^8\dz\theta^{12}\dz\theta^{13} + \\&&
    120\ \sqrt{3}\ \theta^1\dz\theta^4\dz\theta^5\dz\theta^7\dz\theta^8 + 
    36\ \sqrt{3}\ \theta^1\dz\theta^4\dz\theta^5\dz\theta^{11}\dz\theta^{13} + 
    12\ \sqrt{3}\ \theta^1\dz\theta^4\dz\theta^5\dz\theta^{12}\dz\theta^{14} -\\&& 
    720\ \theta^1\dz\theta^4\dz\theta^6\dz\theta^9\dz\theta^{10} - 
    18\ \theta^1\dz\theta^4\dz\theta^6\dz\theta^{13}\dz\theta^{14} - 
    72\ \theta^1\dz\theta^4\dz\theta^8\dz\theta^9\dz\theta^{14} + \\&&
    144\ \theta^1\dz\theta^4\dz\theta^8\dz\theta^{10}\dz\theta^{13} - 
    720\ \sqrt{3}\ \theta^1\dz\theta^5\dz\theta^6\dz\theta^9\dz\theta^{10} - 
    108\ \sqrt{3}\ \theta^1\dz\theta^5\dz\theta^6\dz\theta^{11}\dz\theta^{12} -\\&& 
    6\ \sqrt{3}\ \theta^1\dz\theta^5\dz\theta^6\dz\theta^{13}\dz\theta^{14} - 
    144\ \sqrt{3}\ \theta^1\dz\theta^5\dz\theta^7\dz\theta^9\dz\theta^{12} + 
    192\ \sqrt{3}\ \theta^1\dz\theta^5\dz\theta^7\dz\theta^{10}\dz\theta^{11} - \\&&
    24\ \sqrt{3}\ \theta^1\dz\theta^5\dz\theta^8\dz\theta^9\dz\theta^{14} + 
    48\ \sqrt{3}\ \theta^1\dz\theta^5\dz\theta^8\dz\theta^{10}\dz\theta^{13} + 
    144\ \theta^1\dz\theta^6\dz\theta^7\dz\theta^9\dz\theta^{13} + \\&&
    72\ \theta^1\dz\theta^6\dz\theta^7\dz\theta^{10}\dz\theta^{14} + 
    288\ \theta^1\dz\theta^6\dz\theta^8\dz\theta^9\dz\theta^{11} + 
    216\ \theta^1\dz\theta^6\dz\theta^8\dz\theta^{10}\dz\theta^{12} - \\&&
    54\ \theta^1\dz\theta^9\dz\theta^{10}\dz\theta^{11}\dz\theta^{14} + 
    72\ \theta^1\dz\theta^9\dz\theta^{10}\dz\theta^{12}\dz\theta^{13} - 
    360\ \theta^2\dz\theta^3\dz\theta^4\dz\theta^7\dz\theta^8 - \\&&
    108\ \theta^2\dz\theta^3\dz\theta^4\dz\theta^{11}\dz\theta^{13} - 
    36\ \theta^2\dz\theta^3\dz\theta^4\dz\theta^{12}\dz\theta^{14} + 
    120\ \sqrt{3}\ \theta^2\dz\theta^3\dz\theta^5\dz\theta^7\dz\theta^8 +\\&& 
    36\ \sqrt{3}\ \theta^2\dz\theta^3\dz\theta^5\dz\theta^{11}\dz\theta^{13} + 
    12\ \sqrt{3}\ \theta^2\dz\theta^3\dz\theta^5\dz\theta^{12}\dz\theta^{14} - 
    360\ \theta^2\dz\theta^3\dz\theta^6\dz\theta^9\dz\theta^{10} - \\&&
    162\ \theta^2\dz\theta^3\dz\theta^6\dz\theta^{11}\dz\theta^{12} - 
    18\ \theta^2\dz\theta^3\dz\theta^6\dz\theta^{13}\dz\theta^{14} - 
    216\ \theta^2\dz\theta^3\dz\theta^7\dz\theta^9\dz\theta^{12} + \\&&
    288\ \theta^2\dz\theta^3\dz\theta^7\dz\theta^{10}\dz\theta^{11} - 
    72\ \theta^2\dz\theta^3\dz\theta^8\dz\theta^9\dz\theta^{14} + 
    144\ \theta^2\dz\theta^3\dz\theta^8\dz\theta^{10}\dz\theta^{13} + \\&&
    120\ \sqrt{3}\ \theta^2\dz\theta^4\dz\theta^5\dz\theta^6\dz\theta^8 + 
    48\ \sqrt{3}\ \theta^2\dz\theta^4\dz\theta^5\dz\theta^9\dz\theta^{13} + 
    24\ \sqrt{3}\ \theta^2\dz\theta^4\dz\theta^5\dz\theta^{10}\dz\theta^{14} +\\&& 
    324\ \theta^2\dz\theta^4\dz\theta^7\dz\theta^{11}\dz\theta^{12} + 
    18\ \theta^2\dz\theta^4\dz\theta^7\dz\theta^{13}\dz\theta^{14} + 
    54\ \theta^2\dz\theta^4\dz\theta^8\dz\theta^{11}\dz\theta^{14} - \\&&
    72\ \theta^2\dz\theta^4\dz\theta^8\dz\theta^{12}\dz\theta^{13} - 
    144\ \sqrt{3}\ \theta^2\dz\theta^5\dz\theta^6\dz\theta^9\dz\theta^{12} + 
    192\ \sqrt{3}\ \theta^2\dz\theta^5\dz\theta^6\dz\theta^{10}\dz\theta^{11} -\\&& 
    240\ \sqrt{3}\ \theta^2\dz\theta^5\dz\theta^7\dz\theta^9\dz\theta^{10} - 
    324\ \sqrt{3}\ \theta^2\dz\theta^5\dz\theta^7\dz\theta^{11}\dz\theta^{12} - 
    6\ \sqrt{3}\ \theta^2\dz\theta^5\dz\theta^7\dz\theta^{13}\dz\theta^{14} - \\&&
    18\ \sqrt{3}\ \theta^2\dz\theta^5\dz\theta^8\dz\theta^{11}\dz\theta^{14} + 
    24\ \sqrt{3}\ \theta^2\dz\theta^5\dz\theta^8\dz\theta^{12}\dz\theta^{13} + 
    108\ \theta^2\dz\theta^6\dz\theta^7\dz\theta^{11}\dz\theta^{13} + \\&&
    36\ \theta^2\dz\theta^6\dz\theta^7\dz\theta^{12}\dz\theta^{14} + 
    288\ \theta^2\dz\theta^7\dz\theta^8\dz\theta^9\dz\theta^{11} + 
    216\ \theta^2\dz\theta^7\dz\theta^8\dz\theta^{10}\dz\theta^{12} + \\&&
    24\ \theta^2\dz\theta^9\dz\theta^{11}\dz\theta^{12}\dz\theta^{14} - 
    48\ \theta^2\dz\theta^{10}\dz\theta^{11}\dz\theta^{12}\dz\theta^{13} + 
    120\ \sqrt{3}\ \theta^3\dz\theta^4\dz\theta^5\dz\theta^6\dz\theta^7 + \\&&
    96\ \sqrt{3}\ \theta^3\dz\theta^4\dz\theta^5\dz\theta^9\dz\theta^{11} + 
    72\ \sqrt{3}\ \theta^3\dz\theta^4\dz\theta^5\dz\theta^{10}\dz\theta^{12} + 
    72\ \theta^3\dz\theta^4\dz\theta^6\dz\theta^9\dz\theta^{14} - \\&&
    144\ \theta^3\dz\theta^4\dz\theta^6\dz\theta^{10}\dz\theta^{13} + 
    54\ \theta^3\dz\theta^4\dz\theta^7\dz\theta^{11}\dz\theta^{14} - 
    72\ \theta^3\dz\theta^4\dz\theta^7\dz\theta^{12}\dz\theta^{13} + \\&&
    360\ \theta^3\dz\theta^4\dz\theta^8\dz\theta^9\dz\theta^{10} + 
    162\ \theta^3\dz\theta^4\dz\theta^8\dz\theta^{11}\dz\theta^{12} + 
    72\ \theta^3\dz\theta^4\dz\theta^8\dz\theta^{13}\dz\theta^{14} + \\&&
    24\ \sqrt{3}\ \theta^3\dz\theta^5\dz\theta^6\dz\theta^9\dz\theta^{14} - 
    48\ \sqrt{3}\ \theta^3\dz\theta^5\dz\theta^6\dz\theta^{10}\dz\theta^{13} - 
    18\ \sqrt{3}\ \theta^3\dz\theta^5\dz\theta^7\dz\theta^{11}\dz\theta^{14} +\\&& 
    24\ \sqrt{3}\ \theta^3\dz\theta^5\dz\theta^7\dz\theta^{12}\dz\theta^{13} + 
    120\ \sqrt{3}\ \theta^3\dz\theta^5\dz\theta^8\dz\theta^9\dz\theta^{10} - 
    54\ \sqrt{3}\ \theta^3\dz\theta^5\dz\theta^8\dz\theta^{11}\dz\theta^{12} + \\&&
    108\ \theta^3\dz\theta^6\dz\theta^8\dz\theta^{11}\dz\theta^{13} + 
    36\ \theta^3\dz\theta^6\dz\theta^8\dz\theta^{12}\dz\theta^{14} + 
    144\ \theta^3\dz\theta^7\dz\theta^8\dz\theta^9\dz\theta^{13} + \\&&
    72\ \theta^3\dz\theta^7\dz\theta^8\dz\theta^{10}\dz\theta^{14} - 
    6\ \theta^3\dz\theta^9\dz\theta^{12}\dz\theta^{13}\dz\theta^{14} + 
    12\ \theta^3\dz\theta^{10}\dz\theta^{11}\dz\theta^{13}\dz\theta^{14} -\\&& 
    18\ \sqrt{3}\ \theta^4\dz\theta^5\dz\theta^6\dz\theta^{11}\dz\theta^{14} + 
    24\ \sqrt{3}\ \theta^4\dz\theta^5\dz\theta^6\dz\theta^{12}\dz\theta^{13} - 
    24\ \sqrt{3}\ \theta^4\dz\theta^5\dz\theta^7\dz\theta^9\dz\theta^{14} + \\&&
    48\ \sqrt{3}\ \theta^4\dz\theta^5\dz\theta^7\dz\theta^{10}\dz\theta^{13} - 
    72\ \sqrt{3}\ \theta^4\dz\theta^5\dz\theta^8\dz\theta^9\dz\theta^{12} + 
    96\ \sqrt{3}\ \theta^4\dz\theta^5\dz\theta^8\dz\theta^{10}\dz\theta^{11} + \\&&
    144\ \theta^4\dz\theta^6\dz\theta^8\dz\theta^9\dz\theta^{13} + 
    72\ \theta^4\dz\theta^6\dz\theta^8\dz\theta^{10}\dz\theta^{14} - 
    108\ \theta^4\dz\theta^7\dz\theta^8\dz\theta^{11}\dz\theta^{13} - \\&&
    36\ \theta^4\dz\theta^7\dz\theta^8\dz\theta^{12}\dz\theta^{14} - 
    18\ \theta^4\dz\theta^9\dz\theta^{10}\dz\theta^{13}\dz\theta^{14} + 
    3\ \theta^4\dz\theta^{11}\dz\theta^{12}\dz\theta^{13}\dz\theta^{14} - \\&&
    192\ \sqrt{3}\ \theta^5\dz\theta^6\dz\theta^7\dz\theta^9\dz\theta^{11} - 
    144\ \sqrt{3}\ \theta^5\dz\theta^6\dz\theta^7\dz\theta^{10}\dz\theta^{12} + 
    48\ \sqrt{3}\ \theta^5\dz\theta^6\dz\theta^8\dz\theta^9\dz\theta^{13} + \\&&
    24\ \sqrt{3}\ \theta^5\dz\theta^6\dz\theta^8\dz\theta^{10}\dz\theta^{14} + 
    36\ \sqrt{3}\ \theta^5\dz\theta^7\dz\theta^8\dz\theta^{11}\dz\theta^{13} + 
    12\ \sqrt{3}\ \theta^5\dz\theta^7\dz\theta^8\dz\theta^{12}\dz\theta^{14} + \\&&
    108\ \sqrt{3}\ \theta^5\dz\theta^9\dz\theta^{10}\dz\theta^{11}\dz\theta^{12} - 
    6\ \sqrt{3}\ \theta^5\dz\theta^9\dz\theta^{10}\dz\theta^{13}\dz\theta^
          14 - \sqrt{3}\ \theta^5\dz\theta^{11}\dz\theta^{12}\dz\theta^{13}\dz\theta^{14} - \\&&
    360\ \theta^6\dz\theta^7\dz\theta^8\dz\theta^9\dz\theta^{10} - 
    162\ \theta^6\dz\theta^7\dz\theta^8\dz\theta^{11}\dz\theta^{12} - 
    18\ \theta^6\dz\theta^7\dz\theta^8\dz\theta^{13}\dz\theta^{14} -\\&& 
    108\ \theta^6\dz\theta^9\dz\theta^{10}\dz\theta^{11}\dz\theta^{13} - 
    36\ \theta^6\dz\theta^9\dz\theta^{10}\dz\theta^{12}\dz\theta^{14} + 
    48\ \theta^7\dz\theta^9\dz\theta^{11}\dz\theta^{12}\dz\theta^{13} + \\&&
    24\ \theta^7\dz\theta^{10}\dz\theta^{11}\dz\theta^{12}\dz\theta^{14} - 
    12\ \theta^8\dz\theta^9\dz\theta^{11}\dz\theta^{13}\dz\theta^{14} - 
    6\ \theta^8\dz\theta^{10}\dz\theta^{12}\dz\theta^{13}\dz\theta^{14}.
\end{eqnarray*}

\end{document}